\documentclass[11pt, oneside, a4paper]{scrartcl}

\usepackage[english]{babel}       		
\usepackage{csquotes}

\usepackage[square,numbers]{natbib}        
\RequirePackage{xspace}
\newcommand*{\eg}{\textit{e.g.}\@\xspace}   
\newcommand*{\ie}{\textit{i.e.}\@\xspace}   
\newcommand*{\ale}{\textit{a.e.}\@\xspace}  

\usepackage{graphicx}
\usepackage{amsmath, amssymb, mathrsfs} 	

\newcommand{\diff}{\mathrm{d}}		

\usepackage{scalerel}				
\DeclareMathOperator*{\bigcross}{\scalerel*{\times}{\textstyle\sum}} 
\usepackage{bm}							
\usepackage{upgreek}					
\DeclareMathOperator{\bomega}{\bm{\upomega}}	
\newcommand*{\Su}{\mathbb{S}}					
\newcommand*{\E}{\mathbb{E}}					
\newcommand*{\Var}{\text{Var}}					
\newcommand{\C}{\mathbb{C}}		
\newcommand{\N}{\mathbb{N}}
\newcommand{\R}{\mathbb{R}}

\renewcommand{\phi}{\varphi}		
\renewcommand{\theta}{\vartheta}
\renewcommand{\epsilon}{\varepsilon}

\usepackage{url}
\usepackage{hyperref}
\hypersetup{
	linktocpage,
	colorlinks,
	citecolor=black,  
	filecolor=black,
	linkcolor=black,
	urlcolor=black,   
}

\usepackage{authblk} 		

\usepackage{xcolor}
\definecolor{celltab}{HTML}{EDAB7F}%
\definecolor{mycolor1}{RGB}{0,109,219}%
\definecolor{mycolor2}{RGB}{0,73,73}%
\definecolor{mycolor3}{RGB}{146,0,   0}%
\newcommand{\sae}[1]{\textcolor{mycolor1}{#1}}
\newcommand{\mssae}[1]{\textcolor{mycolor2}{#1}}
\newcommand{\mnssae}[1]{\textcolor{mycolor3}{#1}}

\usepackage[font={small,it}, labelfont=bf]{caption}[2004/07/16] 
\usepackage{subcaption}
\RequirePackage{xparse, tikz}
\RequirePackage{pgfplots}
\pgfplotsset{compat=1.12}
\usetikzlibrary{positioning,calc,shapes,arrows}
\pgfplotsset{every tick label/.append style={font=\scriptsize}}
\newenvironment{customlegend}[1][]{%
	\begingroup
	\csname pgfplots@init@cleared@structures\endcsname
	\pgfplotsset{#1}%
}{%
\csname pgfplots@createlegend\endcsname 
\endgroup
}%
\def\addlegendimage{\csname pgfplots@addlegendimage\endcsname}
\pgfplotsset{
	cycle list={%
		{draw=black,mark=star,solid},
		{draw=black, mark=square,solid},
		{draw=black,mark=+,solid},
		{black,mark=o},
		{draw=black, mark=none,solid}}}

\usepackage{booktabs}  		
\usepackage{multirow}							
\usepackage{colortbl}  							

\usepackage{amsthm}
\theoremstyle{plain}
\newtheorem{theorem}{Theorem}[] 

\theoremstyle{definition}

\newtheorem{example}{Example}

\theoremstyle{remark}

\newtheorem*{remark}{\textbf{Remark}}

\title{\LARGE{Polynomial (chaos) approximation of maximum eigenvalue functions: efficiency and limitations}}
\date{}
\author[1]{Luca Fenzi \thanks{\href{mailto:luca.fenzi@cs.kuleuven.be}{Luca.Fenzi@cs.kuleuven.be}.}}
\author[1]{Wim Michiels \thanks{\href{mailto:wim.michiels@cs.kuleuven.be}{Wim.Michiels@cs.kuleuven.be}.}}
\affil[1]{{\footnotesize Department of Computer Science, NUMA Section, KU Leuven, Belgium}}

\begin{document}

\maketitle

\vspace{-0.1\textwidth}
\begin{abstract}
This paper is concerned with polynomial approximations of the spectral abscissa function (the supremum of the real parts of the eigenvalues) of a parameterized eigenvalue problem, which are closely related to polynomial chaos approximations if the parameters correspond to  realizations of  random variables.
Unlike in existing works, we highlight the major role of the smoothness properties  of the spectral abscissa function. Even if the matrices of the eigenvalue problem are analytic functions of the parameters, the spectral abscissa function may not be everywhere differentiable, even not everywhere Lipschitz continuous, which is related to  multiple rightmost eigenvalues or rightmost eigenvalues with multiplicity higher than one.
The presented analysis demonstrates that the smoothness properties heavily affect the approximation errors of the Galerkin and collocation based polynomial approximations, and the numerical errors of the evaluation of coefficients with integration methods. A documentation of the experiments, conducted on the benchmark problems through the software Chebfun, is publicly available.

\end{abstract}

\section{Introduction}\label{sec:intro}
Polynomial approximation represents a pillar of approximation theory and it is strongly connected with the Polynomial Chaos (PC) method in uncertainty quantification (see,~e.g.~\cite{Trefethen2012} and \cite{LeMaitre2010,Xiu2010,Marelli2017}, respectively, and references therein). The PC method is often used to approximate a quantity of interest as a function of uncertain parameters and to extract relevant statistical information from the approximation.
A topic which is gaining attention, both from a theoretical and application point of view, concerns the PC expansion of eigenvalue functions such as the spectral abscissa function. We refer to, e.g.,~\cite{Ghanem2007, Rahman2011, Sarrouy2012,elman} for the standard eigenvalue problem and \cite{Vermiglio2017} for eigenvalue problems associated with delay differential equations. 

The accuracy of  the polynomial approximation crucially relies on the smoothness property of the function to be approximated. The aim of this work is to show,  in a systematic way, how the different behaviors of rightmost (or dominant) eigenvalues, and in particular spectral abscissa function,  may affect both the quality of polynomial (or PC) based approximations as well as the numerical computation of coefficients. The important link with smoothness properties is barely exploited in the literature.

More precisely, we consider a class of eigenvalue problems inferred from a linear autonomous system of delay differential equations,
\begin{equation}
	\left(\lambda I_n - \sum_{i=0}^h A_i(\omega)e^{-\lambda\tau_i(\omega)}\right)v=0, \ \lambda\in\C,\ v\in\C^n\setminus\{0\},
	\label{eq:eigproblem}
\end{equation}
where $I_n$ is the identity matrix of dimension $n$, $\omega\in\Su\subset\R^D$ models parameters subject to uncertainty, and for every $i\in\{0,\ldots,h\}$, $A_i:\Su\mapsto\R^{n\times n}$ and $\tau_i:\Su\mapsto\R_{\geq0}$ are smooth functions. The eigenvalue problem (\ref{eq:eigproblem}) plays an important role in the stability analysis.  For example, for given $\omega$, the time delay system associated to \eqref{eq:eigproblem},
\begin{equation}
\dot x(t)=\sum_{i=0}^h A_i(\omega)x(t-\tau_i(\omega)),
\end{equation}
 is asymptotically stable  if and only if the real part of the rightmost eigenvalue is negative, or equivalently, if and only if the spectral abscissa
\begin{equation}
	\alpha(\omega):=\max_{\lambda\in\C} \left\{\Re(\lambda): \det\left(\lambda I_n - \sum_{i=0}^h A_i(\omega)e^{-\lambda\tau_i(\omega)}\right)=0 \right\}
	\label{eq:spabs}
\end{equation}
is negative. 

Let $\{p_i(\omega)\}_{i=0}^\infty$ be a degree graded polynomial basis orthogonal w.r.t.~a smooth non-negative function $\rho(\omega)$ defined and normalized on the compact support $\Su\subset\R^D$, \ie $\int_{\Su}\rho(\omega)\diff\omega=1$. 

The spectral abscissa function 
\begin{equation}
	\alpha: \Su \subset\R^D\to \R,\  \omega \mapsto \alpha(\omega),
	\label{eq:spabsfun}
\end{equation}
presents a polynomial expansion w.r.t. the basis $\{p_i(\omega)\}_{i=0}^\infty$
\begin{equation}
	\alpha(\omega)=\sum_{i=0}^\infty c_ip_i(\omega),
	\label{eq:spabsapp}
\end{equation}
where the coefficients $c_i$ are evaluated by
\begin{equation}
	{c}_i=\frac{\langle \alpha, p_i\rangle_{\rho}}{\langle p_i, p_i\rangle_{\rho}},\quad  i\in\N.
	\label{eq:ci}
\end{equation}
The $\rho$-inner product, used in the above equation, is defined for all $f$, $g:\Su\to\R$ such that 
\begin{equation}
	\langle f, g\rangle_{\rho}=\int_{\Su}f(\omega)g(\omega)\rho(\omega)\diff \omega.
\end{equation}
 and the induced $\rho$-norm is determined  by $\|f\|_\rho=\sqrt{\langle f, f\rangle_{\rho}}$.

If $\omega$ is considered as a realization of a real random vector $\bomega$, with  probability density function $\rho(\omega)$, then $\alpha(\bomega)$ is also a random variable and its polynomial expansion
\begin{equation}
	\alpha(\bomega)=\sum_{i=0}^\infty c_ip_i(\bomega),
	\label{eq:spabspce}
\end{equation}
corresponds to a PC expansion of the spectral abscissa function with germ $\bomega$.

The PC expansion permits to compactly define a general random variable, through the chaos coefficients $c_i$, and the germ $\bomega$, which specifies also the ${\rho}$-orthonormal polynomial basis $\{p_i\}_i$. In this context, the chaos coefficients are related to the variance-based sensitivity analysis, as stated in the following theorem.\footnote{Result~\textit{\ref{th:sobol}} in Theorem~\ref{th:pce} is derived by \cite{Crestaux2009}.}

\begin{theorem}\label{th:pce}  Given PC expansion \eqref{eq:spabspce}, the following formulas hold:
	\begin{enumerate}
		\item \label{th:mean} Mean: $\E(\alpha(\bomega))=c_0$.
		\item \label{th:var} Variance: $\Var(\alpha(\bomega))=\sum_{i=1}^\infty c_i^2\langle p_i,p_i\rangle_{\rho}$.
		\item  \label{th:sobol} If the polynomial basis  $\{p_i(\omega)\}_{i=0}^\infty$ is constructed by tensor product of univariate polynomial bases $\{p_{i,d}(\omega_d)\}_{d=1,\ldots,D}^{i=0,\ldots,\infty}$ through the bijective $D$-tupling function $\phi_D: \N^D\to\N$, $(i_1,\ldots,i_D)\mapsto i$ in such a way that
		\begin{equation}\label{eq:approxd2}
			p_i(\omega)=\prod_{d=1}^D p_{i_d,d}(\omega_d),\quad (i_1,\ldots,i_D)=\phi_D^{-1}(i).
		\end{equation}
		Then, the Sobol's sensitivity index $\mathcal{S}_{\mathfrak{h}}(\alpha(\bomega))$ associated to the subset of uncertain parameter $\bomega_{\mathfrak{h}}$, where $\mathfrak{h}=(h_1,\ldots,h_\ell)$ is an element of  the power set of $\{1,\ldots,D\}$ except the empty set,  is determined by:
		\begin{equation}
			\mathcal{S}_{\mathfrak{h}}(\alpha(\bomega))=\sum_{i\in Q_{\mathfrak{h}}}\frac{c_{i}^2\langle p_{i}, p_{i}\rangle_{\rho}}{\Var(\alpha(\bomega))},
		\end{equation}
		where $Q_{\mathfrak{h}}=\{i\in\N\ :\  \phi_D^{-1}(i)=(i_{h_1},\ldots,i_{h_D}), \text{ and } (i_{h_1},\ldots,i_{h_\ell})>0,\  (i_{h_{\ell+1}},\ldots,i_{h_D})=0 \}$.
		For $d\in\{0,\ldots,D\}$, the effect of the $d$th parameter $\omega_d$ can be quantified by the total order sensitivity index
		\begin{equation}
			\mathcal{S}_{d}(\alpha(\bomega))=\sum_{\substack{\mathfrak{h}=(h_1,\ldots,h_\ell),\\ \exists 1\leq j\leq \ell \ h_j=d}}\mathcal{S}_{\mathfrak{h}}(\alpha(\bomega)).
		\end{equation}
	\end{enumerate}
\end{theorem}

The present work highlights the fundamental role of the behavior of the spectral abscissa function $\alpha(\omega)$ in \eqref{eq:spabsfun} w.r.t.~truncations of its polynomial expansion \eqref{eq:spabsapp}  and the corresponding PC expansion \eqref{eq:spabspce} of $\alpha(\bomega)$.

\medskip

The paper is organized as follows. In Section~\ref{sec:behavior}, we briefly review the behaviors of spectral abscissa functions furnishing text examples which are going to be analyzed in Section~\ref{sec:d1}, where the polynomial approximation for parameter eigenvalue problems with $D=1$  is investigated. Then, in Section~\ref{sec:d2}, the polynomial approximation for $D>1$ is analyzed through some numerical experiments.  Finally, we end with some concluding remarks in Section~\ref{sec:conclusions}. A MATLAB tutorial with the numerical experiments here presented can be found in \cite{Fenzi2018}.

\section{Smoothness property of the spectral abscissa function}\label{sec:behavior}
Due to the link between PC expansion  \eqref{eq:spabspce} of $\alpha(\bomega)$ and the polynomial approximation of spectral abscissa function \eqref{eq:spabsfun}, the smoothness properties of this function play a fundamental role in the analysis.
In this section, we characterize the different behaviors of the spectral abscissa furnishing benchmark examples, which are analyzed in the following  Section~\ref{sec:d1}. The reader is referred to \cite{Fenzi2017, Michiels2017}  for relevant results on the behavior of spectral abscissa functions. 

The different behaviors of the spectral abscissa can be inferred by the analysis of the active eigenvalues, which are defined as eigenvalues ($\lambda$ in \eqref{eq:eigproblem}) with real part equal to the spectral abscissa ($\Re(\lambda)=\alpha$) and non-negative imaginary part ($\Im(\lambda)\geq0$). 

The spectral abscissa is a continuous function which is smooth \ale in $\Su$, however in a set of measure zero, it may not be differentiable due to the presence of more than one active eigenvalue (counted with multiplicity).

\begin{description}
	\item[\sae{SAE}] - \sae{Simple Active Eigenvalue}. For all $\omega\in\Su$, there is only one active eigenvalue, whose algebraic and geometric multiplicity is equal to $1$. The spectral abscissa is smooth over all the uncertain domain $\Su$, \ie $\alpha\in\mathcal{C}^\infty(\Su,\R)$. 
	\item[\mssae{MSSAEs}] - \mssae{Multiple (Semi)-Simple Active Eigenvalues}. For some $\omega\in\Su$, there is more than one active eigenvalue. If the multiplicity of an active eigenvalue is greater than one, then it is semi-simple (\ie the algebraic multiplicity is equal to the geometric one). In this case, the spectral abscissa may not be everywhere differentiable, however it is Lipschitz continuous in $\Su$.
	\item[\mnssae{MNSSAEs}] - \mnssae{Multiple Non-Semi-Simple Active Eigenvalues}. There exist $\omega\in\Su$, where the algebraic multiplicity of an active eigenvalue is greater than the geometric one (\ie it is  non-semi-simple). In this case, the spectral abscissa is, in general, not even everywhere locally Lipschitz continuous.  
\end{description}

\begin{example}\label{ex:d1} The characterization of the behaviors of the spectral abscissa function can be illustrated by the following parameter eigenvalue problems, where  $\omega\in[-1,1]$,
	\begin{align}
		\tag{\text{\sae{SAE}}}\label{eq:SAE1}	\left(\lambda I_2-\begin{pmatrix}
			e^\omega & 0 \\ 0 & -1
		\end{pmatrix}\right)v&=0;\\ 
		\tag{\mssae{MSSAEs}}\label{eq:MSSAE1}	\left(\lambda I_2-\begin{pmatrix}
			\omega & 0 \\ 0 & 0
		\end{pmatrix}\right)v&=0;\\
		\tag{\mnssae{MNSSAEs}}\label{eq:MNSSAE1}	\left(\lambda I_2-\begin{pmatrix}
			0 & \omega \\ 1 & 0
		\end{pmatrix}\right)v&=0.
	\end{align}

	The real part of the spectra are shown in Figure~\ref{fig:spectrum}. The first eigenvalue problem, \eqref{eq:SAE1}, presents a smooth spectral abscissa function, the second one, \eqref{eq:MSSAE1}, presents a piecewise linear function (for $\omega=0$ the active eigenvalue is multiple semi-simple), while the last one, \eqref{eq:MNSSAE1}, shows the splitting of a double eigenvalue at $\omega=0$ with a square root behavior for $\omega\in[0,1]$. 

	\begin{figure}[h]
		\begin{subfigure}[h]{0.32\linewidth}
			\centering
			\begin{tikzpicture}

\begin{axis}[%
width=0.8\textwidth,
height=0.4\textwidth,
at={(0\textwidth,0\textwidth)},
scale only axis,
xmin=-1,
xmax=1,
xlabel={$\omega$},
ymin=-1.5,
ymax=3,
ytick={-1,0,...,2,3},
ylabel={$\Re(\lambda)$},
title style={font=\bfseries, yshift=-0.2cm},
axis background/.style={fill=white},
title={\eqref{eq:SAE1}},
 ylabel style = {yshift=-0.3cm, font=\footnotesize},
 xlabel style = {yshift=+0.1cm, font=\footnotesize},
]
\addplot [color=mycolor1,solid,line width=2.0pt,forget plot]
  table[row sep=crcr]{%
-1	0.367879441171442\\
-0.97979797979798	0.375386926934554\\
-0.95959595959596	0.383047621429048\\
-0.939393939393939	0.390864651256307\\
-0.919191919191919	0.398841206823714\\
-0.898989898989899	0.406980543646768\\
-0.878787878787879	0.415285983677781\\
-0.858585858585859	0.423760916661676\\
-0.838383838383838	0.432408801519472\\
-0.818181818181818	0.441233167759984\\
-0.797979797979798	0.450237616920347\\
-0.777777777777778	0.459425824035927\\
-0.757575757575758	0.468801539140235\\
-0.737373737373737	0.478368588795451\\
-0.717171717171717	0.488130877654176\\
-0.696969696969697	0.498092390053061\\
-0.676767676767677	0.508257191638956\\
-0.656565656565657	0.518629431028247\\
-0.636363636363636	0.52921334150005\\
-0.616161616161616	0.540013242723965\\
-0.595959595959596	0.551033542523085\\
-0.575757575757576	0.562278738672988\\
-0.555555555555556	0.573753420737433\\
-0.535353535353535	0.585462271941531\\
-0.515151515151515	0.59741007108313\\
-0.494949494949495	0.609601694483216\\
-0.474747474747475	0.62204211797611\\
-0.454545454545455	0.634736418940282\\
-0.434343434343434	0.647689778370613\\
-0.414141414141414	0.660907482992939\\
-0.393939393939394	0.674394927421755\\
-0.373737373737374	0.688157616361946\\
-0.353535353535353	0.702201166855453\\
-0.333333333333333	0.716531310573789\\
-0.313131313131313	0.731153896157336\\
-0.292929292929293	0.746074891602382\\
-0.272727272727273	0.761300386696874\\
-0.252525252525252	0.776836595505876\\
-0.232323232323232	0.792689858907749\\
-0.212121212121212	0.80886664718209\\
-0.191919191919192	0.82537356265048\\
-0.171717171717172	0.84221734237113\\
-0.151515151515151	0.859404860888509\\
-0.131313131313131	0.876943133039094\\
-0.111111111111111	0.89483931681437\\
-0.0909090909090909	0.913100716282262\\
-0.0707070707070707	0.931734784568187\\
-0.0505050505050505	0.950749126896934\\
-0.0303030303030303	0.97015150369663\\
-0.0101010101010101	0.989949833766045\\
0.0101010101010102	1.01015219750654\\
0.0303030303030303	1.03076684021994\\
0.0505050505050506	1.05180217547379\\
0.0707070707070707	1.07326678853516\\
0.0909090909090908	1.09516943987466\\
0.111111111111111	1.11751906874186\\
0.131313131313131	1.14032479681373\\
0.151515151515152	1.16359593191751\\
0.171717171717172	1.18734197182957\\
0.191919191919192	1.21157260815182\\
0.212121212121212	1.23629773026713\\
0.232323232323232	1.2615274293756\\
0.252525252525253	1.28727200261311\\
0.272727272727273	1.31354195725395\\
0.292929292929293	1.34034801499921\\
0.313131313131313	1.36770111635268\\
0.333333333333333	1.39561242508609\\
0.353535353535354	1.4240933327954\\
0.373737373737374	1.45315546355014\\
0.393939393939394	1.48281067863759\\
0.414141414141414	1.51307108140382\\
0.434343434343434	1.54394902219345\\
0.454545454545455	1.57545710339032\\
0.474747474747475	1.60760818456091\\
0.494949494949495	1.64041538770285\\
0.515151515151515	1.67389210260041\\
0.535353535353535	1.70805199228938\\
0.555555555555556	1.74290899863346\\
0.575757575757576	1.77847734801437\\
0.595959595959596	1.8147715571382\\
0.616161616161616	1.85180643896017\\
0.636363636363636	1.88959710873031\\
0.656565656565657	1.92815899016254\\
0.676767676767677	1.96750782172967\\
0.696969696969697	2.00765966308675\\
0.717171717171717	2.04863090162566\\
0.737373737373737	2.09043825916337\\
0.757575757575758	2.13309879876667\\
0.777777777777778	2.17662993171625\\
0.797979797979798	2.22104942461286\\
0.818181818181818	2.26637540662847\\
0.838383838383838	2.31262637690544\\
0.858585858585859	2.35982121210669\\
0.878787878787879	2.40797917411992\\
0.898989898989899	2.45711991791906\\
0.919191919191919	2.50726349958618\\
0.939393939393939	2.5584303844971\\
0.95959595959596	2.61064145567402\\
0.97979797979798	2.6639180223086\\
1	2.71828182845905\\
};
\addplot [color=black,dashed,line width=1.2pt,forget plot]
  table[row sep=crcr]{%
-1	0.367879441171442\\
-0.97979797979798	0.375386926934554\\
-0.95959595959596	0.383047621429048\\
-0.939393939393939	0.390864651256307\\
-0.919191919191919	0.398841206823714\\
-0.898989898989899	0.406980543646768\\
-0.878787878787879	0.415285983677781\\
-0.858585858585859	0.423760916661676\\
-0.838383838383838	0.432408801519472\\
-0.818181818181818	0.441233167759984\\
-0.797979797979798	0.450237616920347\\
-0.777777777777778	0.459425824035927\\
-0.757575757575758	0.468801539140235\\
-0.737373737373737	0.478368588795451\\
-0.717171717171717	0.488130877654176\\
-0.696969696969697	0.498092390053061\\
-0.676767676767677	0.508257191638956\\
-0.656565656565657	0.518629431028247\\
-0.636363636363636	0.52921334150005\\
-0.616161616161616	0.540013242723965\\
-0.595959595959596	0.551033542523085\\
-0.575757575757576	0.562278738672988\\
-0.555555555555556	0.573753420737433\\
-0.535353535353535	0.585462271941531\\
-0.515151515151515	0.59741007108313\\
-0.494949494949495	0.609601694483216\\
-0.474747474747475	0.62204211797611\\
-0.454545454545455	0.634736418940282\\
-0.434343434343434	0.647689778370613\\
-0.414141414141414	0.660907482992939\\
-0.393939393939394	0.674394927421755\\
-0.373737373737374	0.688157616361946\\
-0.353535353535353	0.702201166855453\\
-0.333333333333333	0.716531310573789\\
-0.313131313131313	0.731153896157336\\
-0.292929292929293	0.746074891602382\\
-0.272727272727273	0.761300386696874\\
-0.252525252525252	0.776836595505876\\
-0.232323232323232	0.792689858907749\\
-0.212121212121212	0.80886664718209\\
-0.191919191919192	0.82537356265048\\
-0.171717171717172	0.84221734237113\\
-0.151515151515151	0.859404860888509\\
-0.131313131313131	0.876943133039094\\
-0.111111111111111	0.89483931681437\\
-0.0909090909090909	0.913100716282262\\
-0.0707070707070707	0.931734784568187\\
-0.0505050505050505	0.950749126896934\\
-0.0303030303030303	0.97015150369663\\
-0.0101010101010101	0.989949833766045\\
0.0101010101010102	1.01015219750654\\
0.0303030303030303	1.03076684021994\\
0.0505050505050506	1.05180217547379\\
0.0707070707070707	1.07326678853516\\
0.0909090909090908	1.09516943987466\\
0.111111111111111	1.11751906874186\\
0.131313131313131	1.14032479681373\\
0.151515151515152	1.16359593191751\\
0.171717171717172	1.18734197182957\\
0.191919191919192	1.21157260815182\\
0.212121212121212	1.23629773026713\\
0.232323232323232	1.2615274293756\\
0.252525252525253	1.28727200261311\\
0.272727272727273	1.31354195725395\\
0.292929292929293	1.34034801499921\\
0.313131313131313	1.36770111635268\\
0.333333333333333	1.39561242508609\\
0.353535353535354	1.4240933327954\\
0.373737373737374	1.45315546355014\\
0.393939393939394	1.48281067863759\\
0.414141414141414	1.51307108140382\\
0.434343434343434	1.54394902219345\\
0.454545454545455	1.57545710339032\\
0.474747474747475	1.60760818456091\\
0.494949494949495	1.64041538770285\\
0.515151515151515	1.67389210260041\\
0.535353535353535	1.70805199228938\\
0.555555555555556	1.74290899863346\\
0.575757575757576	1.77847734801437\\
0.595959595959596	1.8147715571382\\
0.616161616161616	1.85180643896017\\
0.636363636363636	1.88959710873031\\
0.656565656565657	1.92815899016254\\
0.676767676767677	1.96750782172967\\
0.696969696969697	2.00765966308675\\
0.717171717171717	2.04863090162566\\
0.737373737373737	2.09043825916337\\
0.757575757575758	2.13309879876667\\
0.777777777777778	2.17662993171625\\
0.797979797979798	2.22104942461286\\
0.818181818181818	2.26637540662847\\
0.838383838383838	2.31262637690544\\
0.858585858585859	2.35982121210669\\
0.878787878787879	2.40797917411992\\
0.898989898989899	2.45711991791906\\
0.919191919191919	2.50726349958618\\
0.939393939393939	2.5584303844971\\
0.95959595959596	2.61064145567402\\
0.97979797979798	2.6639180223086\\
1	2.71828182845905\\
};
\addplot [color=white!30!black,dashed,line width=1.2pt,forget plot]
  table[row sep=crcr]{%
-1	-1\\
1	-1\\
};
\end{axis}
\end{tikzpicture}%
		\end{subfigure}\quad
		\begin{subfigure}[h]{0.32\linewidth}
			\centering
			\begin{tikzpicture}

\begin{axis}[%
width=0.8\textwidth,
height=0.4\textwidth,
at={(0\textwidth,0\textwidth)},
scale only axis,
xmin=-1,
xmax=1,
xlabel={$\omega$},
ymin=-1,
ymax=1,
ylabel={$\Re(\lambda)$},
axis background/.style={fill=white},
title style={font=\bfseries, yshift=-0.2cm},
title={\eqref{eq:MSSAE1}},
ylabel style = {yshift=-0.3cm, font=\footnotesize},
xlabel style = {yshift=+0.1cm,  font=\footnotesize}
]
\addplot [color=mycolor2,solid,line width=2.0pt,forget plot]
  table[row sep=crcr]{%
-1	0\\
0	0\\
1	1\\
};
\addplot [color=black,dashed,line width=1.2pt,forget plot]
  table[row sep=crcr]{%
-1	-1\\
1	1\\
};
\addplot [color=white!30!black,dashed,line width=1.2pt,forget plot]
  table[row sep=crcr]{%
-1	0\\
1	0\\
};
\end{axis}
\end{tikzpicture}%
		\end{subfigure}	\quad
		\begin{subfigure}[h]{0.32\linewidth}
			\centering
			\begin{tikzpicture}

\begin{axis}[%
width=0.8\textwidth,
height=0.4\textwidth,
at={(0\textwidth,0\textwidth)},
scale only axis,
xmin=-1,
xmax=1,
xlabel={$\omega$},
ymin=-1,
ymax=1,
ylabel={$\Re(\lambda)$},
axis background/.style={fill=white},
title style={font=\bfseries, yshift=-0.2cm},
title={\eqref{eq:MNSSAE1}},
ylabel style = {yshift=-0.3cm, font=\footnotesize},
xlabel style = {yshift=+0.1cm, font=\footnotesize}
]
\addplot [color=mycolor3,solid,line width=2.0pt,forget plot]
  table[row sep=crcr]{%
-1	0\\
0	0\\
0.01	0.1\\
0.02	0.14142135623731\\
0.03	0.173205080756888\\
0.04	0.2\\
0.05	0.223606797749979\\
0.06	0.244948974278318\\
0.07	0.264575131106459\\
0.08	0.282842712474619\\
0.09	0.3\\
0.1	0.316227766016838\\
0.11	0.33166247903554\\
0.12	0.346410161513775\\
0.13	0.360555127546399\\
0.14	0.374165738677394\\
0.15	0.387298334620742\\
0.16	0.4\\
0.17	0.412310562561766\\
0.18	0.424264068711929\\
0.19	0.435889894354067\\
0.2	0.447213595499958\\
0.21	0.458257569495584\\
0.22	0.469041575982343\\
0.23	0.479583152331272\\
0.24	0.489897948556636\\
0.25	0.5\\
0.26	0.509901951359279\\
0.27	0.519615242270663\\
0.28	0.529150262212918\\
0.29	0.53851648071345\\
0.3	0.547722557505166\\
0.31	0.556776436283002\\
0.32	0.565685424949238\\
0.33	0.574456264653803\\
0.34	0.58309518948453\\
0.35	0.591607978309962\\
0.36	0.6\\
0.37	0.608276253029822\\
0.38	0.616441400296898\\
0.39	0.62449979983984\\
0.4	0.632455532033676\\
0.41	0.640312423743285\\
0.42	0.648074069840786\\
0.43	0.6557438524302\\
0.44	0.66332495807108\\
0.45	0.670820393249937\\
0.46	0.678232998312527\\
0.47	0.685565460040104\\
0.48	0.692820323027551\\
0.49	0.7\\
0.5	0.707106781186548\\
0.51	0.714142842854285\\
0.52	0.721110255092798\\
0.53	0.728010988928052\\
0.54	0.734846922834953\\
0.55	0.741619848709566\\
0.56	0.748331477354788\\
0.57	0.754983443527075\\
0.58	0.761577310586391\\
0.59	0.768114574786861\\
0.6	0.774596669241483\\
0.61	0.781024967590665\\
0.62	0.787400787401181\\
0.63	0.793725393319377\\
0.64	0.8\\
0.65	0.806225774829855\\
0.66	0.812403840463596\\
0.67	0.818535277187245\\
0.68	0.824621125123532\\
0.69	0.830662386291807\\
0.7	0.836660026534076\\
0.71	0.842614977317636\\
0.72	0.848528137423857\\
0.73	0.854400374531753\\
0.74	0.860232526704263\\
0.75	0.866025403784439\\
0.76	0.871779788708135\\
0.77	0.877496438739212\\
0.78	0.883176086632785\\
0.79	0.888819441731559\\
0.8	0.894427190999916\\
0.81	0.9\\
0.82	0.905538513813742\\
0.83	0.91104335791443\\
0.84	0.916515138991168\\
0.85	0.921954445729289\\
0.86	0.92736184954957\\
0.87	0.932737905308881\\
0.88	0.938083151964686\\
0.89	0.94339811320566\\
0.9	0.948683298050514\\
0.91	0.953939201416946\\
0.92	0.959166304662544\\
0.93	0.964365076099295\\
0.94	0.969535971483266\\
0.95	0.974679434480896\\
0.96	0.979795897113271\\
0.97	0.98488578017961\\
0.98	0.989949493661167\\
0.99	0.99498743710662\\
1	1\\
};
\addplot [color=black,dashed,line width=1.2pt,forget plot]
  table[row sep=crcr]{%
-1	0\\
0	0\\
};
\addplot [color=black,dashed,line width=1.2pt,forget plot]
  table[row sep=crcr]{%
0	0\\
0.01	0.1\\
0.02	0.14142135623731\\
0.03	0.173205080756888\\
0.04	0.2\\
0.05	0.223606797749979\\
0.06	0.244948974278318\\
0.07	0.264575131106459\\
0.08	0.282842712474619\\
0.09	0.3\\
0.1	0.316227766016838\\
0.11	0.33166247903554\\
0.12	0.346410161513775\\
0.13	0.360555127546399\\
0.14	0.374165738677394\\
0.15	0.387298334620742\\
0.16	0.4\\
0.17	0.412310562561766\\
0.18	0.424264068711929\\
0.19	0.435889894354067\\
0.2	0.447213595499958\\
0.21	0.458257569495584\\
0.22	0.469041575982343\\
0.23	0.479583152331272\\
0.24	0.489897948556636\\
0.25	0.5\\
0.26	0.509901951359279\\
0.27	0.519615242270663\\
0.28	0.529150262212918\\
0.29	0.53851648071345\\
0.3	0.547722557505166\\
0.31	0.556776436283002\\
0.32	0.565685424949238\\
0.33	0.574456264653803\\
0.34	0.58309518948453\\
0.35	0.591607978309962\\
0.36	0.6\\
0.37	0.608276253029822\\
0.38	0.616441400296898\\
0.39	0.62449979983984\\
0.4	0.632455532033676\\
0.41	0.640312423743285\\
0.42	0.648074069840786\\
0.43	0.6557438524302\\
0.44	0.66332495807108\\
0.45	0.670820393249937\\
0.46	0.678232998312527\\
0.47	0.685565460040104\\
0.48	0.692820323027551\\
0.49	0.7\\
0.5	0.707106781186548\\
0.51	0.714142842854285\\
0.52	0.721110255092798\\
0.53	0.728010988928052\\
0.54	0.734846922834953\\
0.55	0.741619848709566\\
0.56	0.748331477354788\\
0.57	0.754983443527075\\
0.58	0.761577310586391\\
0.59	0.768114574786861\\
0.6	0.774596669241483\\
0.61	0.781024967590665\\
0.62	0.787400787401181\\
0.63	0.793725393319377\\
0.64	0.8\\
0.65	0.806225774829855\\
0.66	0.812403840463596\\
0.67	0.818535277187245\\
0.68	0.824621125123532\\
0.69	0.830662386291807\\
0.7	0.836660026534076\\
0.71	0.842614977317636\\
0.72	0.848528137423857\\
0.73	0.854400374531753\\
0.74	0.860232526704263\\
0.75	0.866025403784439\\
0.76	0.871779788708135\\
0.77	0.877496438739212\\
0.78	0.883176086632785\\
0.79	0.888819441731559\\
0.8	0.894427190999916\\
0.81	0.9\\
0.82	0.905538513813742\\
0.83	0.91104335791443\\
0.84	0.916515138991168\\
0.85	0.921954445729289\\
0.86	0.92736184954957\\
0.87	0.932737905308881\\
0.88	0.938083151964686\\
0.89	0.94339811320566\\
0.9	0.948683298050514\\
0.91	0.953939201416946\\
0.92	0.959166304662544\\
0.93	0.964365076099295\\
0.94	0.969535971483266\\
0.95	0.974679434480896\\
0.96	0.979795897113271\\
0.97	0.98488578017961\\
0.98	0.989949493661167\\
0.99	0.99498743710662\\
1	1\\
};
\addplot [color=white!30!black,dashed,line width=1.2pt,forget plot]
  table[row sep=crcr]{%
0	-0\\
0.01	-0.1\\
0.02	-0.14142135623731\\
0.03	-0.173205080756888\\
0.04	-0.2\\
0.05	-0.223606797749979\\
0.06	-0.244948974278318\\
0.07	-0.264575131106459\\
0.08	-0.282842712474619\\
0.09	-0.3\\
0.1	-0.316227766016838\\
0.11	-0.33166247903554\\
0.12	-0.346410161513775\\
0.13	-0.360555127546399\\
0.14	-0.374165738677394\\
0.15	-0.387298334620742\\
0.16	-0.4\\
0.17	-0.412310562561766\\
0.18	-0.424264068711929\\
0.19	-0.435889894354067\\
0.2	-0.447213595499958\\
0.21	-0.458257569495584\\
0.22	-0.469041575982343\\
0.23	-0.479583152331272\\
0.24	-0.489897948556636\\
0.25	-0.5\\
0.26	-0.509901951359279\\
0.27	-0.519615242270663\\
0.28	-0.529150262212918\\
0.29	-0.53851648071345\\
0.3	-0.547722557505166\\
0.31	-0.556776436283002\\
0.32	-0.565685424949238\\
0.33	-0.574456264653803\\
0.34	-0.58309518948453\\
0.35	-0.591607978309962\\
0.36	-0.6\\
0.37	-0.608276253029822\\
0.38	-0.616441400296898\\
0.39	-0.62449979983984\\
0.4	-0.632455532033676\\
0.41	-0.640312423743285\\
0.42	-0.648074069840786\\
0.43	-0.6557438524302\\
0.44	-0.66332495807108\\
0.45	-0.670820393249937\\
0.46	-0.678232998312527\\
0.47	-0.685565460040104\\
0.48	-0.692820323027551\\
0.49	-0.7\\
0.5	-0.707106781186548\\
0.51	-0.714142842854285\\
0.52	-0.721110255092798\\
0.53	-0.728010988928052\\
0.54	-0.734846922834953\\
0.55	-0.741619848709566\\
0.56	-0.748331477354788\\
0.57	-0.754983443527075\\
0.58	-0.761577310586391\\
0.59	-0.768114574786861\\
0.6	-0.774596669241483\\
0.61	-0.781024967590665\\
0.62	-0.787400787401181\\
0.63	-0.793725393319377\\
0.64	-0.8\\
0.65	-0.806225774829855\\
0.66	-0.812403840463596\\
0.67	-0.818535277187245\\
0.68	-0.824621125123532\\
0.69	-0.830662386291807\\
0.7	-0.836660026534076\\
0.71	-0.842614977317636\\
0.72	-0.848528137423857\\
0.73	-0.854400374531753\\
0.74	-0.860232526704263\\
0.75	-0.866025403784439\\
0.76	-0.871779788708135\\
0.77	-0.877496438739212\\
0.78	-0.883176086632785\\
0.79	-0.888819441731559\\
0.8	-0.894427190999916\\
0.81	-0.9\\
0.82	-0.905538513813742\\
0.83	-0.91104335791443\\
0.84	-0.916515138991168\\
0.85	-0.921954445729289\\
0.86	-0.92736184954957\\
0.87	-0.932737905308881\\
0.88	-0.938083151964686\\
0.89	-0.94339811320566\\
0.9	-0.948683298050514\\
0.91	-0.953939201416946\\
0.92	-0.959166304662544\\
0.93	-0.964365076099295\\
0.94	-0.969535971483266\\
0.95	-0.974679434480896\\
0.96	-0.979795897113271\\
0.97	-0.98488578017961\\
0.98	-0.989949493661167\\
0.99	-0.99498743710662\\
1	-1\\
};
\end{axis}
\end{tikzpicture}%
		\end{subfigure}
		\vspace{-0.05\textwidth}
		\caption{Real part of the spectra varying the $\omega$ parameter in $[-1,1]$ of the eigenvalue problems defined in Example~\ref{ex:d1}. The spectral abscissa functions are highlighted.} 
		\label{fig:spectrum}
		
	\end{figure}
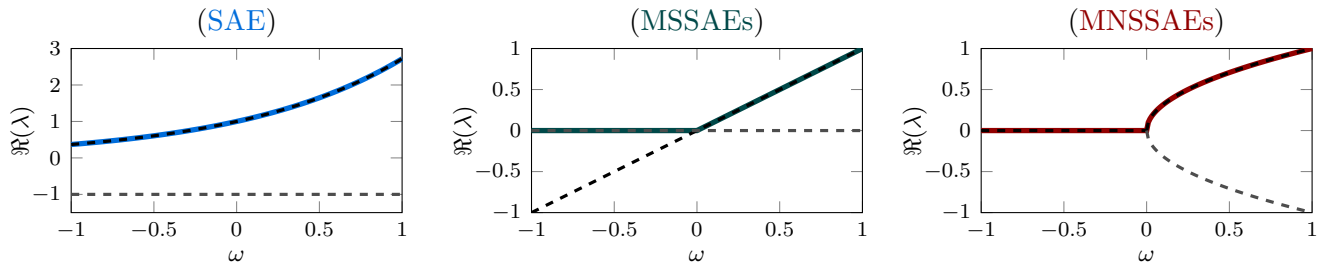

\end{example}

\section{Analysis of the  polynomial approximation for $D=1$} \label{sec:d1}
In practice, for computational reasons, polynomial expansion \eqref{eq:spabsapp} is approximated by finitely many terms 
\begin{equation}
	\alpha_P(\omega)=\sum_{i=0}^P \tilde{c}_ip_i(\omega).
	\label{eq:spabsapprox}
\end{equation} 
Several methods exist to compute the coefficients $\tilde{c}_i$; here we focus on the Galerkin and collocation approaches, analyzed in Sections~\ref{sec:Galerkin} and \ref{sec:collocation}, respectively.

In this section, the  $L^\infty$ convergence rate of polynomial approximation \eqref{eq:spabsapprox} is analyzed for the model problems in Example~\ref{ex:d1}. (We consider $\rho(\omega)=1/2$ uniformly distributed in $\Su=[-1,1]$, and therefore the polynomial basis orthogonal to the $\rho$-inner product is determined by Legendre polynomials.)

\subsection{Galerkin approach}\label{sec:Galerkin}
Given a finite polynomial basis $\{p_i\}_{i=0}^P$, the Galerkin approach finds an approximation \eqref{eq:spabsapprox} of spectral abscissa function \eqref{eq:spabsfun} such that the residual is orthogonal w.r.t.~the polynomial basis, in formula
\begin{equation}\label{eq:galerkin}	
	\langle \alpha-\alpha_P,p_i\rangle_\rho=0, \quad i=0,\ldots,P.
\end{equation} 
This leads to $\tilde{c}_i=c_i$, where the coefficients $c_i$ are defined in \eqref{eq:ci}. Hence, the Galerkin approach is nothing else than a truncation up to order $P$ of polynomial series \eqref{eq:spabsapp}. Moreover, it provides the optimal approximation in the $\rho$-norm, \ie 
\begin{equation}
	\|\alpha-\alpha_P\|_\rho=\sqrt{\int_\Su \left(\alpha(\omega)-\alpha_P(\omega)\right)^2\rho(\omega)\diff\omega}
\end{equation} 
is minimized. 

This approximation of $\alpha(\omega)$ corresponds to the stochastic Galerkin approximation of $\alpha(\bomega)$ in the PC theory.

The analyses of the convergence consider, first of all, the truncation error up to order $P$, and then, the numerical error due to the computation of the coefficients $c_i$.  

\subsubsection{Approximation error}
Assuming that the coefficients $c_i$ are correctly evaluated, the error bounds of polynomial approximation \eqref{eq:spabsapprox} obtained by the Galerkin approach with Legendre polynomials $\{p_i\}_{i=0}^P$ (with $\rho(\omega)=1/2$ for $\omega\in[-1,1]$) are determined by the following theorem, which follows from \cite{Wang2011}.
\begin{theorem}\label{th:convgalerkin}
	Let $\alpha$, $\alpha^\prime$, $\ldots$ $\alpha^{(k-1)}$ be absolutely continuous on $\Su=[-1,1]$ for $k>1$ and \begin{equation}\label{eq:boundedvariation1}
		\int_{-1}^{1} \frac{\lvert \alpha^{(k+1)}(\omega) \rvert}{\sqrt{1-\omega^2}}\diff\omega=V_k<\infty.
	\end{equation} Then, for each $P>k+1$ and $\omega\in\Su$, the following relation holds
	\begin{equation}\label{eq:inequality2}
		\| \alpha -\alpha_P\|_{\infty}\leq\frac{V_k\sqrt{\pi/2}}{(k-1)(P-k)^{k-0.5}}.\end{equation}
\end{theorem}

Condition \eqref{eq:boundedvariation1} corresponds to the fact that the $k$th derivative, $\alpha^{(k)}$, presents bounded variation w.r.t.~the Chebyshev weight function. 

The level of smoothness, indicated by the maximum $k$ which satisfies the assumption of Theorem~\ref{th:convgalerkin}, heavily determines the convergence behavior of polynomial approximation \eqref{eq:spabsapprox}. Only the \sae{SAE} case verifies the assumptions of Theorem~\ref{th:convgalerkin}, and it satisfies them for all $k\in\N$, ensuring a convergence rate faster than $\mathcal{O}(P^{-k})$, for all $k\in\N$.
In fact, it is shown in \cite{Wang2011}  that the convergence rate is at least geometric, based on analytic extension of the real function on the complex plane; and such extension always exists for a  real analytic function.

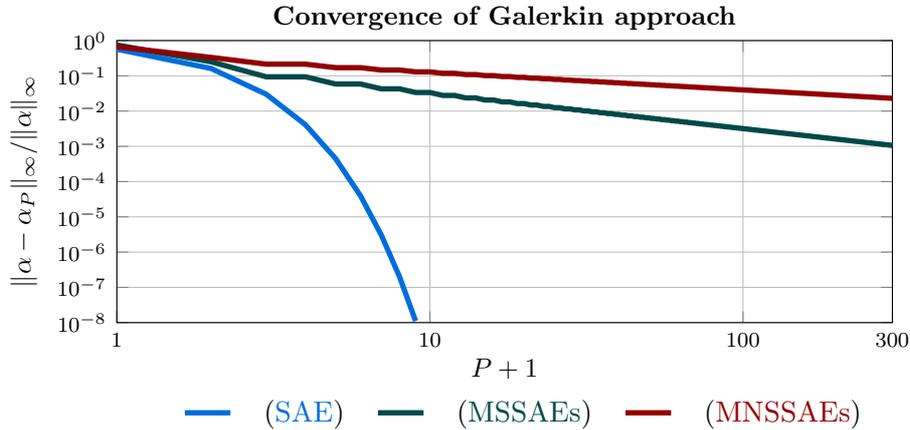
\begin{figure}[!h]
	\begin{center}
		\include{GalD1}
		\vspace{-0.03\textwidth}  
		\begin{tikzpicture}
		\hspace{0.05\linewidth}
		\begin{customlegend}[legend columns=3,legend style={align=left,draw=none,column sep=2ex},legend entries={\eqref{eq:SAE1},\eqref{eq:MSSAE1},\eqref{eq:MNSSAE1}}]
		\addlegendimage{color=mycolor1, mark=none,solid,line width=2.0pt,line legend}
		\addlegendimage{color=mycolor2, mark=none,solid,line width=2.0pt}   
		\addlegendimage{color=mycolor3, mark=none,solid,line width=2.0pt}
		\end{customlegend}
		\end{tikzpicture}
		\caption{Relative error  of polynomial approximation \eqref{eq:spabsapprox} obtained by the Galerkin approach up to order~$P$ of the spectral abscissa associated to the parameter eigenvalue problems  of Example~\ref{ex:d1}.} 
		\label{fig:d1galerkin}
	\end{center}
\end{figure}

Considering Example~\ref{ex:d1}, the truncation error obtained by the Galerkin approach for the spectral abscissa of parameter eigenvalue problem \eqref{eq:SAE1} presents a spectral convergence, $\mathcal{O}(P^{-P})$, as shown in Figure~\ref{fig:d1galerkin}.  The convergence rates of \eqref{eq:MSSAE1} and \eqref{eq:MNSSAE1} are approximately $\mathcal{O}(P^{-1})$, and $\mathcal{O}(P^{-0.5})$, respectively. The coefficients of the polynomial approximation are analytically evaluated in Appendix~\ref{app:ex1}.

\subsubsection{Numerical error} \label{sec:gd1num}

It is not always possible to analytically compute the coefficients $c_i=\tilde{c}_i$ with \eqref{eq:ci}, therefore an integration method based on $M+1$ nodes can be used to approximate the integrals $\langle \alpha, p_i\rangle_\rho$ and consequently the coefficients, denoted by $\tilde{c}_i^M$. This type of approach is known as Non Intrusive Spectral Projection in the PC framework. 

In this section, we approximate the coefficients of the Galerkin polynomial approximation $\alpha_P^M$ with  the following integration methods based on $M+1$ points:
\begin{description}
	\item[Classical integration methods] based on an equally spaced discretization of $\Su$, which considers the extremes of $\Su$, \ie $-1$ and $1$.  In particular, we focus on extended trapezoidal rule  and extended Simpson's rule. For the latter method, $M$ is required to be an even number.
	\item[Interpolatory quadrature rules] which approximates an integral by integrating the interpolant of its integrand, where the degree of the polynomial interpolant is at most $M$. We consider Clenshaw-Curtis and Gauss quadrature rules based on Chebyshev and Legendre points, respectively. (Further information on the interpolatory properties of Chebyshev points are given in the upcoming Section~\ref{sec:collocation}.)
\end{description}

For smooth function, including the \sae{SAE} case, the extended trapezoidal and Simpson's rules provide an error of the order $\mathcal{O}(M^{-2})$ and $\mathcal{O}(M^{-4})$, respectively. 

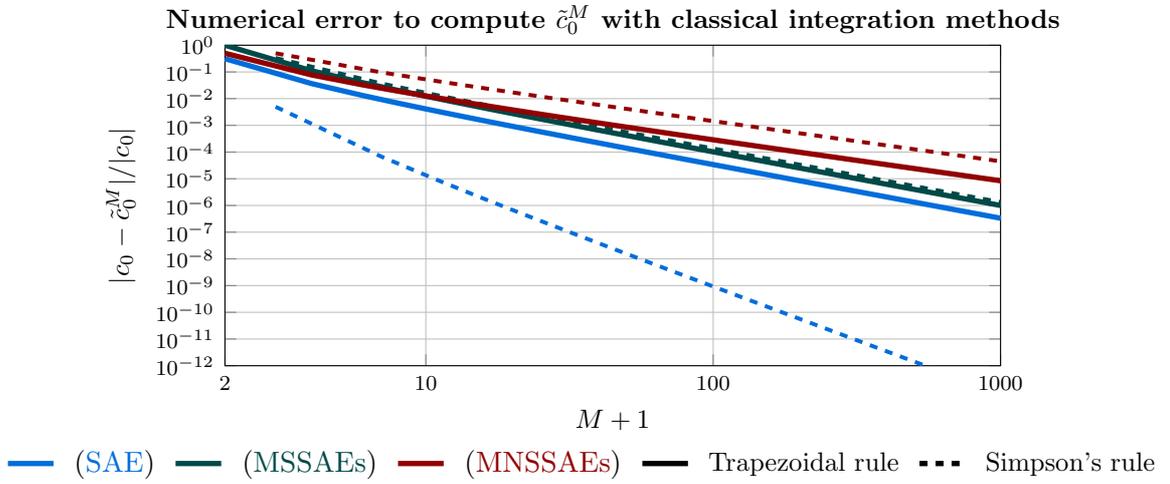
\begin{figure}[!h]
	\begin{center}
		\include{c0D1_trapz}
		\vspace{-0.03\textwidth}  
		\begin{tikzpicture}
		\begin{customlegend}[legend columns=5,legend style={align=left,draw=none,column sep=1ex},legend entries={\eqref{eq:SAE1},\eqref{eq:MSSAE1},\eqref{eq:MNSSAE1},\small{Trapezoidal rule},  \small{Simpson's rule}}]
		\addlegendimage{color=mycolor1, mark=none,solid,line width=2.0pt,line legend}
		\addlegendimage{color=mycolor2, mark=none,solid,line width=2.0pt}   
		\addlegendimage{color=mycolor3, mark=none,solid,line width=2.0pt}
		\addlegendimage{color=black, mark=none,solid,line width=2.0pt}
		\addlegendimage{color=black, mark=none,dashed,line width=2.0pt}
		\end{customlegend}
		\end{tikzpicture}
		\caption{Error to compute an approximation of the first coefficient $\tilde{c}_0^M$ of \eqref{eq:spabsapprox} with classical integration methods (extended Trapezoidal and Simpson's rules) based on $M+1$ equally spaced points for problems of Example~\ref{ex:d1}. Only the slowest convergence rates, for the different methods and problems, are displayed.}
		\label{fig:c0d1_trapz}
	\end{center}
\end{figure}

 Figure~\ref{fig:c0d1_trapz} shows the numerical errors of classical integration methods, to evaluate the first coefficient $c_0$ of polynomial expansion \eqref{eq:spabsapp} for Example~\ref{ex:d1}. The convergences for \eqref{eq:SAE1} follow the theoretical error bounds. For  \eqref{eq:MSSAE1} and \eqref{eq:MNSSAE1}, the convergence rates for both classical integration methods are approximately $\mathcal{O}(M^{-2})$ and $\mathcal{O}(M^{-1.5})$, respectively.
  
Whenever the integrand is smooth enough to be well-approximated by a polynomial, the interpolatory quadrature rules perform better than the classical integration rule; as stated in the following theorem, which combines results from \cite{Trefethen2012, Xiang2012, Xiang2016}. 
\begin{theorem}\label{th:integrationd1} For an integer $k\geq1$, let $\alpha$, $\alpha^\prime$, $\ldots$ $\alpha^{(k-1)}$ be absolutely continuous on $\Su=[-1,1]$ and let $(\alpha(\omega) p_i(\omega))^{(k)}$ present bounded total variation $V_i$ with $i\in\N$. Then, the Clenshaw-Curtis quadrature rule, for the approximation of  $c_i$ with $M+1$ points, \ie  $\tilde{c}_i^M$ satisfies
	\begin{equation} \label{eq:ccint}
	\lvert c_i-\tilde{c}_i^M\rvert\leq\frac{32}{15}\frac{V_i}{\pi k (M-k)^k}, \text{ for }M>k,
	\end{equation}
	and  $(M+1)$ Gauss quadrature satisfies
	\begin{equation} \label{eq:errata}
	\lvert c_i-\tilde{c}_i^M\rvert\leq\frac{32}{15}\frac{V_i}{\pi k (2M+1-k)^{k}}, \text{ for }M>2k+1.
	\end{equation}
	
Moreover, for $k=1$, the Clenshaw-Curtis quadrature rule presents a convergence rate of $\mathcal{O}(M^{-2})$, while the Gauss quadrature error is at most of size   $\mathcal{O}(M^{-2}\ln(M))$. 
	
\end{theorem}

The theorem can be applied to  \sae{SAE}, where it provides a convergence rate faster than $\mathcal{O}(M^{-k})$ for all $k\in\N$, for both interpolatory quadrature rules. Moreover the theorem ensures that in the \mssae{MSSAEs}  Clenshaw-Curtis quadrature has a rate of convergence of  $\mathcal{O}(M^{-2})$, while the Gauss quadrature errors  converges, at least, as $\mathcal{O}(M^{-2}\ln(M))$.

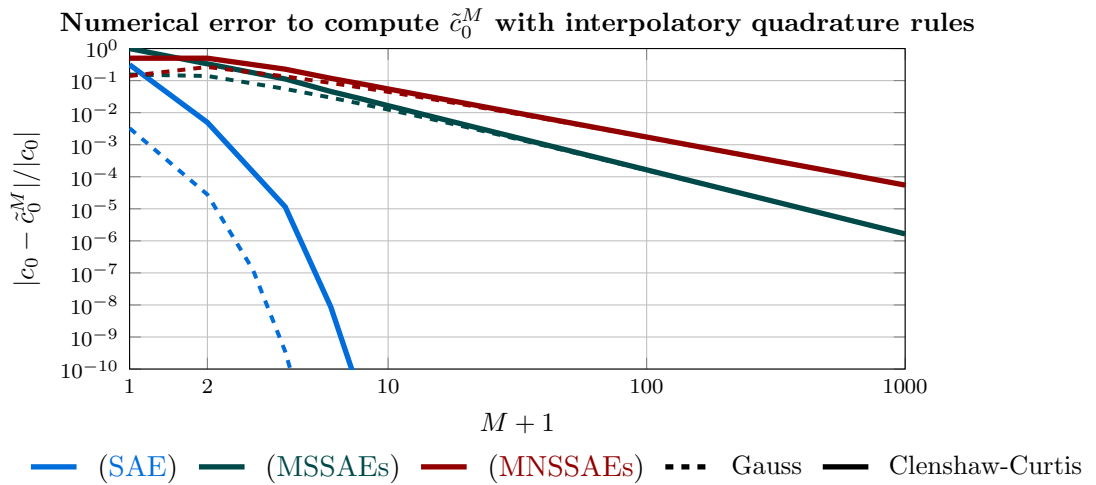
\begin{figure}[!h]
	\begin{center}
		\include{c0D1}
		\vspace{-0.03\textwidth}  
		\begin{tikzpicture}
		\hspace{0.05\linewidth}
		\begin{customlegend}[legend columns=5,legend style={align=left,draw=none,column sep=1ex},legend entries={\eqref{eq:SAE1},\eqref{eq:MSSAE1},\eqref{eq:MNSSAE1}, \small{Gauss}, \small{Clenshaw-Curtis}}]
		\addlegendimage{color=mycolor1, mark=none,solid,line width=2.0pt,line legend}
		\addlegendimage{color=mycolor2, mark=none,solid,line width=2.0pt}   
		\addlegendimage{color=mycolor3, mark=none,solid,line width=2.0pt}
		\addlegendimage{color=black, mark=none,dashed,line width=2.0pt}
		\addlegendimage{color=black, mark=none,solid,line width=2.0pt}
		\end{customlegend}
		\end{tikzpicture}
		\caption{Error to compute an approximation of the first coefficient $\tilde{c}_0^M$ of \eqref{eq:spabsapprox} through Clenshaw-Curtis and Gauss quadrature rules with $(M+1)$ points for problems of Example~\ref{ex:d1}. Only the slowest convergence rates, for the different methods and problems, are displayed.}
		\label{fig:c0d1}
	\end{center}
\end{figure}

Figure~\ref{fig:c0d1} illustrates the numerical error induced by the approximation with interpolatory quadrature rules of the first coefficient $c_0$ of Example~\ref{ex:d1}.  \eqref{eq:SAE1} converges with an order of $\mathcal{O}(M^{-M})$, improving the convergence rate of the classical integration methods. For the non-smooth cases, \ie \eqref{eq:MSSAE1}, and \eqref{eq:MNSSAE1}, the convergence rates of interpolatory quadrature rules are similar to the ones obtained by classical integration methods (Figure~\ref{fig:c0d1_trapz}).

From this numerical experiment, Clenshaw-Curtis and Gauss quadrature rules present similar convergence rates. Hence, the bound of Theorem~\ref{th:integrationd1} for  the \mssae{MSSAEs} are optimal for the Clenshaw-Curtis quadrature rule, while they are conservative for Gauss quadrature rule (as already observed in  \cite{Xiang2012} and \cite{Xiang2016}).

\begin{remark} 
		 When the coefficients  $\tilde{c}_i^M$ of  $\alpha_P^M$ are approximated by interpolatory quadrature rules, it is advised to set $M\geq P$ as explained with further details in the following Section~\ref{sec:coupling}.
\end{remark}
\begin{remark} 
		 In the PC framework, and in particular for Theorem~\ref{th:pce}, the first coefficient $\tilde{c}_0^M$ corresponds to  approximating the mean  of the spectral abscissas in Example~\ref{ex:d1}, where $\omega$ is a realization of random variable uniformly distributed in $[-1,1]$. The corresponding convergence rates for the mean (Figures~\ref{fig:c0d1_trapz} and \ref{fig:c0d1}) provide only insights on the numerical error of the integration method, and are not meaningful for the truncation error of the polynomial approximation, \ie they do not depend on $P$.
\end{remark}

\subsection{Collocation approach}\label{sec:collocation}
The collocation approach determines the coefficients of \eqref{eq:spabsapprox} by interpolation on $P+1$ points, $\{\zeta_i\}_{i=0}^P$:
\begin{equation}\label{eq:collocation}
	\alpha_P(\zeta_i)=\alpha(\zeta_i), \quad \text{for all }\ i\in\{0,\ldots,P\}.
\end{equation}
The coefficients $\tilde{c}_i$ of the polynomial approximation $\alpha_P$ can be computed  solving a linear system of $P+1$ equations, and this can be often done with a negligible numerical error.

A widely used choice of interpolating nodes for the interval $[-1,1]$ are Chebyshev nodes. The polynomial interpolant  satisfies error bounds, similar to Theorem~\ref{th:convgalerkin}, and provides a near-best approximation in $L^\infty$ sense, as stated in the following theorem, whose proof can be found \eg in \cite{Trefethen2012}. 

\begin{theorem}\label{th:convcollocation}
	If $\alpha$, $\alpha^\prime$, $\ldots$ $\alpha^{(k-1)}$ are absolutely continuous on $\Su=[-1,1]$ and $\alpha^{(k)}$ presents bounded total variation $V$ for some $k\geq1$, then the Chebyshev interpolant $\alpha_P$ of degree $P>k$, , satisfies
	\begin{equation}\label{eq:inequality3}
		\| \alpha-\alpha_P\|_{\infty}\leq\frac{4V}{\pi k(P-k)^{k}}.
	\end{equation}
	Moreover, if $\alpha_P^\star$ is the best polynomial approximation of order less than or equal to $P$, then 
	\begin{equation}\label{eq:nearbest}
		\| \alpha-\alpha_P\|_{\infty}\leq\left(2+\frac{2}{\pi}\log(P+1)\right)\| \alpha-\alpha_P^\star\|_{\infty}.
	\end{equation}
\end{theorem}

Theorem~\ref{th:convcollocation} ensures that the collocation approach on Chebyshev nodes converge, at least, as $\mathcal{O}(P^{-1})$ and faster than $\mathcal{O}(P^{-k})$ for all $k\in\N$, for the \mssae{MSSAEs} and \sae{SAE} cases, respectively.

\begin{figure}[!h]
	\begin{center}
		\include{BestD1}
		\vspace{-0.03\textwidth}  
		\begin{tikzpicture}
		\hspace{0.05\linewidth}
		\begin{customlegend}[legend columns=5,legend style={align=left,draw=none,column sep=1ex},legend entries={\eqref{eq:SAE1},\eqref{eq:MSSAE1},\eqref{eq:MNSSAE1}, best $\alpha_P^\star$, near-best $\alpha_P$}]
		\addlegendimage{color=mycolor1, mark=none,solid,line width=2.0pt,line legend}
		\addlegendimage{color=mycolor2, mark=none,solid,line width=2.0pt}   
		\addlegendimage{color=mycolor3, mark=none,solid,line width=2.0pt}
		\addlegendimage{color=black, mark=none,dotted,line width=2.0pt}   
		\addlegendimage{color=black, mark=none,solid,line width=2.0pt}
		\end{customlegend}
		\end{tikzpicture}
		\caption{Relative error of polynomial approximation \eqref{eq:spabsapprox} obtained by interpolation on $P+1$ Chebyshev points (near-best approximation) for the spectral abscissa associated to the parameter eigenvalue problems of Example~\ref{ex:d1}. For completeness the best polynomial approximation in $L^\infty$ sense is shown. Only the odd values of $P$ for \eqref{eq:MSSAE1} and \eqref{eq:MNSSAE1} are considered, since they furnish a slower convergence rate.} 
		\label{fig:bestnearbestd1}
	\end{center}
\end{figure}
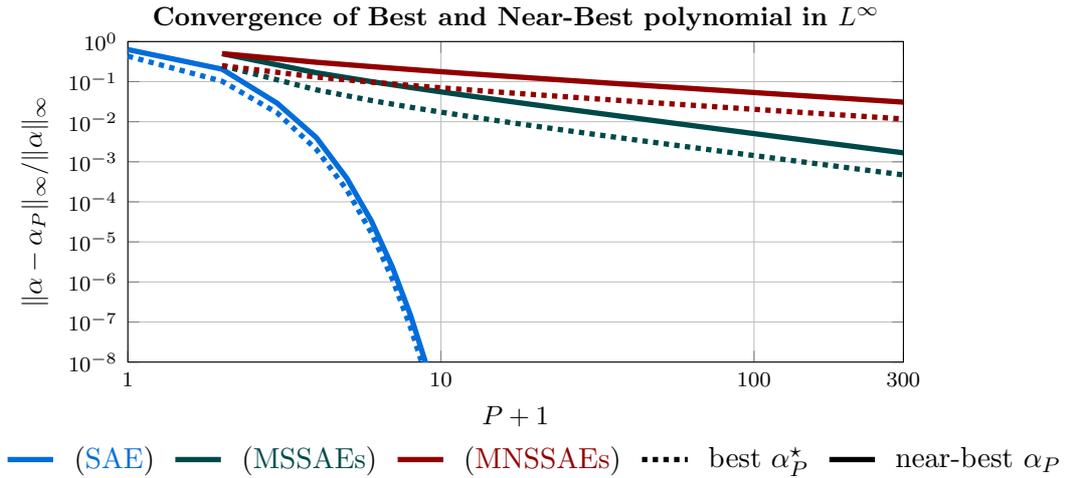 

The computation of the polynomial interpolant on Chebyshev points can be achieved trough Fast Fourier Transform based algorithm, which maps interpolating values onto coefficients of a  polynomial approximation in the Chebyshev polynomial basis. The MATLAB software Chebfun is considered to evaluate the Chebyshev interpolants for the model problems described in Example~\ref{ex:d1}, and can be used to evaluate also the best $L^\infty$ polynomial approximation.  Figure~\ref{fig:bestnearbestd1},  other than indicating the interpolating error (near-best polynomial approximation), shows the convergence rates of the best $L^\infty$ polynomial approximation.
The convergence rates are analogous to the ones obtained by the Galerkin approach (Section~\ref{sec:Galerkin}, Figure~\ref{fig:d1galerkin}) and agree with Theorem~\ref{th:convcollocation} for \sae{SAE} and \mssae{MNSSAEs}.

\begin{remark} If we consider the PC framework, where $\omega$ in Example~\ref{ex:d1} is a realization of $\bomega$, uniformly distributed in $[-1,1]$, and we want to apply Theorem~\ref{th:pce}, a polynomial transformation is needed in order to convert the Chebyshev expansion coefficients into the Legendre coefficients (using  \eg the method in \cite{Townsend2017}). The transformation does not affect the results shown in Figure~\ref{fig:bestnearbestd1}, since the numerical error, also in this case, is negligible. 
\end{remark}

\section{Experiments on the polynomial approximation for $D>1$}\label{sec:d2}
In this section, first of all, we generalize polynomial approximation \eqref{eq:spabsapprox} for $D>1$,  in such a way that the coefficients $\tilde{c}_i$ can be evaluated by the Galerkin and collocation approaches, through  the corresponding formula \eqref{eq:galerkin} and \eqref{eq:collocation}. Then, we consider examples of $D=2$ parameters eigenvalue problems, where the spectral abscissa functions present the different behaviors characterized in Section~\ref{sec:behavior}. The Galerkin and collocation approaches are, hence, applied for the latter benchmark in Sections~\ref{sec:gald2} and \ref{sec:col2}, respectively.

We focus on polynomial approximation using a degree graded polynomial basis $\{p_i(\omega)\}_{i=0}^\infty$, constructed by products of univariate degree graded polynomial bases $\{p_{i,d}(\omega_d)\}_{d=1,\ldots,D}^{i=0,\ldots,\infty}$. The degree grading of the polynomial basis $\{p_i\}_{i}$ follows a norm $\|\cdot\|$ on $\N^D$, which can be associated to a $D$-tupling function $\phi_D$. In this way, the basis is obtained by formula \eqref{eq:approxd2}, and the generalization of polynomial approximation \eqref{eq:spabsapprox} for the multivariate case is given by
\begin{equation}
\alpha_P(\omega)=\sum_{\|\phi_D^{-1}(i)\|\leq P_d,\ i\in\N} \tilde{c}_ip_i(\omega)=\sum_{i=0}^P \tilde{c}_ip_i(\omega),\label{eq:spapproxd2}
\end{equation}
where $P+1$ is equal to the number of polynomials in $\{p_i\}_{i}$ such that the multivariate degree is less than or equal to $P_d\geq0$. 

Two popular choice of norms are the $1$-norm, $\|\phi_D^{-1}(i)\|_1=\sum_{d=1}^D i_d$, and the $\infty$-norm, $\|\phi_D^{-1}(i)\|_\infty=\max_{d=1,\ldots,D} i_d$ , which are associated to the total and maximal degrees, respectively. In these cases, the number of coefficients $\tilde{c}_i$ in \eqref{eq:spapproxd2} satisfies $P+1={P_d+D \choose D}$ and $P+1=(P_d+1)^D$ for the total and maximal degree, respectively.\footnote{The interested reader is referred to \cite{Szudzik2017} and \cite{Trefethen2017} for further  information on  $D$-tupling functions  and on multivariate degree, respectively.}

 If the bases $\{p_{i,d}(\omega_d)\}_{d,i}$  are orthogonal w.r.t.~$\rho_d(\omega_d)$ with $\omega_d\in\Su_d\subset\R$, then the $D$-dimensional polynomial basis $\{p_i(\omega)\}_{i=0}^\infty$ is orthogonal to the normalized weight function
\begin{equation}\label{eq:orthoD}
	\rho(\omega)=\prod_{d=1}^D\rho_d(\omega_d),\quad \omega\in\Su=\bigcross_{d=1}^D\Su_d\subset\R^D.
\end{equation}
In the PC theory, probability density function \eqref{eq:orthoD} corresponds to the assumption that the random vector $\bomega$ is constituted by $D$ independent random variable $\bomega_d$, $d=1,\ldots,D$. In this framework, the multivariate polynomial degree determines the truncation scheme, and in particular, the total degree corresponds to the standard truncation scheme.

In the following Sections~\ref{sec:gald2}  and ~\ref{sec:col2}, we analyze the Galerkin and collocation approaches, respectively, on the following benchmark examples with $D=2$, already studied in \cite{Fenzi2017}. 

\begin{example}\label{ex:d2} We consider the spectral abscissa functions associated to the oscillator with feedback delay  system 	
\begin{equation}\label{eq:spring}
	\ddot{x}(t)=-\omega_1^2x(t)-2\omega_1\omega_2 \dot{x}(t) +K_1x(t-1)+K_2\dot{x}(t-1),
\end{equation}
such that $x(t)$ is the normalized position defined for  time $t\in[-1,\infty)$, $\omega_1\in\Su_1=[0.9,1.1]$ and $\omega_2\in\Su_2=[0.1,0.2]$ are angular frequency and damping ratio, respectively, while $(K_1,K_2)$ describes the control force which acts with a delay of $\tau=1$. The control variables are set equal to the values of Table~\ref{tab:d2}, in this way we illustrate the different behaviors of the spectral abscissa function (Section~\ref{sec:behavior}) for parameter eigenvalue problems with $D=2$. 

\begin{table}[!h]
	\caption{Numerical values of the control parameters $K_1$, $K_2$ for system \eqref{eq:spring} corresponding to different behaviors of the spectral abscissa for $(\omega_1,\omega_2)\in\Su_1\times\Su_2$.\footnotesize{(The table is equal to Table 2 in \cite{Fenzi2017})}}
	\label{tab:d2}
	\begin{center}
		\begin{tabular}{lccc}
			\toprule
			& \sae{SAE} 	& \mssae{MSSAEs} & \mnssae{MNSSAEs}\\
			\midrule
			$K_1$ 	& 0.2 	& 0.5105  & 0.6179\\ 
			$K_2$ 	&0.2 	& -0.0918 & -0.0072\\ 
			\bottomrule
		\end{tabular}
	\end{center}
\end{table}

For this problem, the spectral abscissa functions are not known analytically, even though it is possible to compute their values given $\omega\in\Su=\Su_1\times\Su_2$. In particular, the convergence rates in $L^\infty$ of Figures~\ref{fig:gd2_approx_err},~\ref{fig:countercheck},  and \ref{fig:collocationd2} are  computed w.r.t.~$10^6$ equidistant points in $\Su$.

W.l.g. we embed the linear transformation $[-1,1]^2\to\Su$ in system \eqref{eq:spring}, and in the associated parameter eigenvalue problems. In this way, we do not consider shifting and rescaling in the cubature rules, in the interpolation and in the polynomial bases.  
\end{example}

\begin{figure}[h]
	\begin{center}
		\begin{subfigure}[h]{0.45\linewidth}
			\centering
			\definecolor{Cantor}{rgb}{   0,   0,0}%

\begin{tikzpicture}


  \draw [dotted] (0,0) -- (0,1.5) node [above right=-0.05cm]  {};  
  \draw [dotted] (1.5,0) -- (0,3) node [above right=-0.05cm]  {};  
  \draw [dotted] (3,0) -- (0,4.5) node [above right=-0.05cm]  {}; 
  \draw [dotted] (4.5,0) -- (0.75,5) node [above right=-0.05cm]  {}; 
  \draw [dotted] (5,1.25) -- (2,5) node [above right=-0.05cm]  {}; 
  \draw [dotted] (5,3.3333333333333333333) -- (3,5) node [above right=-0.05cm]  {}; 

\clip (-1,-0.5) rectangle (5.5,5.5); 
\pgftransformcm{1}{0}{0}{1}{\pgfpoint{0cm}{0cm}}
\draw[style=help lines,dashed] (0,0) grid[step=1.5cm] (5,5); 

\foreach \x in {0,1,...,3}{
	\foreach \y in {0,1,...,3}{
		\node[draw,circle,inner sep=1pt,fill, gray] at (1.5*\x,1.5*\y) {};
		\node[below left=-0.05cm, gray]  at (1.5*\x,1.5*\y) {\tiny (\x,\y)};
	}
}

	\node [above right=-0.05cm,-latex,Cantor] at (0,0) {\scriptsize $0$};
    \node [above right=-0.05cm,-latex,Cantor] at (0,1.5) {\scriptsize $1$};
    \draw [->,>=stealth, very thick,Cantor] (0,1.5) -- (1.5,0) node [above right=-0.05cm]  {\scriptsize $2$};
    \node [above right=-0.05cm,-latex,Cantor] at (0,3) {\scriptsize $3$};
    \draw [->,>=stealth, very thick,Cantor]  (0,3)  -- (1.5,1.5) node [above right=-0.05cm]  {\scriptsize $4$};
    \draw [->,>=stealth, very thick,Cantor]  (1.5,1.5)  -- (3,0) node [above right=-0.05cm]  {\scriptsize $5$};
    \node [above right=-0.05cm,-latex,Cantor] at (0,4.5) {\scriptsize $6$};
    \draw [->,>=stealth, very thick,Cantor]  (0,4.5)  -- (1.5,3) node [above right=-0.05cm]  {\scriptsize $7$};
    \draw [->,>=stealth, very thick,Cantor]  (1.5,3)  -- (3,1.5) node [above right=-0.05cm]  {\scriptsize $8$};
    \draw [->,>=stealth, very thick,Cantor]  (3,1.5)  -- (4.5,0) node [above right=-0.05cm]  {\scriptsize $9$};
    \node [above right=-0.05cm,-latex,Cantor] at (1.5,4.5) {\scriptsize $11$};
    \draw [->,>=stealth, very thick,Cantor]  (1.5,4.5)  -- (3,3) node [above right=-0.05cm]  {\scriptsize $12$};
    \draw [->,>=stealth, very thick,Cantor]  (3,3)  -- (4.5,1.5) node [above right=-0.05cm]  {\scriptsize $13$};
    \node [above right=-0.05cm,-latex,Cantor] at (3,4.5) {\scriptsize $17$};
    \draw [->,>=stealth, very thick,Cantor]  (3,4.5)  -- (4.5,3) node [above right=-0.05cm]  {\scriptsize $18$};
    \node [above right=-0.05cm,-latex,Cantor] at (4.5,4.5) {\scriptsize $24$};

    \draw [->,>=stealth, very thick,Cantor]  (1,5)  -- (1.5,4.5) node [above right=-0.05cm]  {};
    \draw [->,>=stealth, very thick,Cantor]  (2.5,5)  -- (3,4.5) node [above right=-0.05cm]  {};
    \draw [->,>=stealth, very thick,Cantor]  (4,5)  -- (4.5,4.5) node [above right=-0.05cm]  {};
    \draw [very thick,Cantor]  (4.5,4.5)  -- (5,4) node [above right=-0.05cm]  {};
    \draw [very thick,Cantor]  (4.5,3)  -- (5,2.5) node [above right=-0.05cm]  {};
    \draw [very thick,Cantor]  (4.5,1.5)  -- (5,1) node [above right=-0.05cm]  {};

  \node[above,font=\small\bfseries] at (2.3,5) {$\pi_1$\textcolor{gray}{$(i_1,i_2)$}$=i$};

  \end{tikzpicture}	
		\end{subfigure}\quad
		\begin{subfigure}[h]{0.45\linewidth}
			\centering
			\definecolor{Cantor}{rgb}{   0,   0,0}%

\begin{tikzpicture}
  \draw [dotted] (0,0) -- (0,1.5) node [above right=-0.05cm]  {};  
  \draw [dotted] (1.5,0) -- (0,3) node [above right=-0.05cm]  {};  
  \draw [dotted] (3,0) -- (0,4.5) node [above right=-0.05cm]  {}; 
  \draw [dotted] (4.5,0) -- (0.75,5) node [above right=-0.05cm]  {}; 
  \draw [dotted] (5,1.25) -- (2,5) node [above right=-0.05cm]  {}; 
  \draw [dotted] (5,3.3333333333333333333) -- (3,5) node [above right=-0.05cm]  {};

\clip (-1,-0.5) rectangle (5.5,5.5); 
\pgftransformcm{1}{0}{0}{1}{\pgfpoint{0cm}{0cm}}
\draw[style=help lines,dashed] (0,0) grid[step=1.5cm] (5,5); 

\foreach \x in {0,1,...,3}{
	\foreach \y in {0,1,...,3}{
		\node[draw,circle,inner sep=1pt,fill, gray] at (1.5*\x,1.5*\y) {};
		\node[below left=-0.05cm, gray]  at (1.5*\x,1.5*\y) {\tiny (\x,\y)};
	}
}

\node [above right=-0.05cm,-latex,Cantor] at (0,0) {\scriptsize $0$};
\node [above right=-0.05cm,-latex,Cantor] at (0,1.5) {\scriptsize $1$};
\draw [->,>=stealth, very thick, Cantor]  (0,1.5)   -- (1.5,1.5) node [above right=-0.05cm]  {\scriptsize $2$};
\draw [->,>=stealth, very thick, Cantor]  (1.5,1.5)   -- (1.5,0) node [above right=-0.05cm]  {\scriptsize $3$};
\node [above right=-0.05cm,-latex,Cantor] at (0,3) {\scriptsize $4$};
\draw [->,>=stealth, very thick,Cantor]  (0,3)   -- (1.5,3) node [above right=-0.05cm]  {\scriptsize $5$};
\draw [->,>=stealth, very thick,Cantor]  (1.5,3)   -- (3,3) node [above right=-0.05cm]  {\scriptsize $6$};
\draw [->,>=stealth, very thick,Cantor]  (3,3)   -- (3,1.5) node [above right=-0.05cm]  {\scriptsize $7$};
\draw [->,>=stealth, very thick,Cantor]  (3,1.5)   -- (3,0) node [above right=-0.05cm]  {\scriptsize $8$};
\node [above right=-0.05cm,-latex,Cantor] at (0,4.5) {\scriptsize $9$};
\draw [->,>=stealth, very thick,Cantor]  (0,4.5)   -- (1.5,4.5) node [above right=-0.05cm]  {\scriptsize $10$};
\draw [->,>=stealth, very thick,Cantor]  (1.5,4.5)   -- (3,4.5) node [above right=-0.05cm]  {\scriptsize $11$};
\draw [->,>=stealth, very thick,Cantor]  (3,4.5)   -- (4.5,4.5) node [above right=-0.05cm]  {\scriptsize $12$};
\draw [->,>=stealth, very thick,Cantor]  (4.5,4.5)   -- (4.5,3) node [above right=-0.05cm]  {\scriptsize $13$};
\draw [->,>=stealth, very thick,Cantor]  (4.5,3)   -- (4.5,1.5) node [above right=-0.05cm]  {\scriptsize $14$};
\draw [->,>=stealth, very thick,Cantor]  (4.5,2)   -- (4.5,0) node [above right=-0.05cm]  {\scriptsize $15$};

  \node[above,font=\small\bfseries] at (2.3,5) {$\pi_\infty$\textcolor{gray}{$(i_1,i_2)$}$=i$};      
  \end{tikzpicture}	
		\end{subfigure}
	\end{center}
	\vspace{-0.05\textwidth}
	\caption{(Left pane) Cantor pairing function \eqref{eq:piCantor}. (Right pane) Rosenberg-Strong pairing function \eqref{eq:piRosenberg}.  } 
	\label{fig:pi}
	
\end{figure}
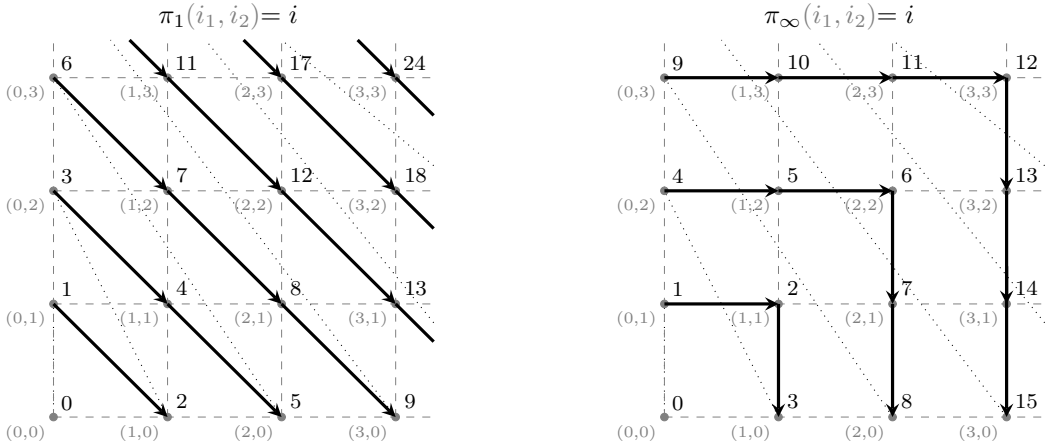

The degree grading of the polynomial basis, (considered in the following sections) for $D=2$ parameter eigenvalue problems of Example~\ref{ex:d2}, are specified by:
\begin{description}
	\item[Total Degree] is determined by the $1$-norm,  such that $\|(i_1,i_2)\|_1=i_1+i_2$. The total degree is associated with a pairing function (\ie a $D=2$-tupling function)  that assigns consecutive numbers to points along diagonals of $\N\times\N$, \eg the Cantor pairing function $\pi_1(i_1,i_2)$ represented in the left pane of Figure~\ref{fig:pi} and defined by
	\begin{equation}\label{eq:piCantor}
	\pi_1(i_1,i_2)=\frac{i_1^2+3i_1+2i_1i_2+i_2+i_2^2}{2}.
	\end{equation}
	The polynomial approximation, obtained from the basis corresponding to the pairing function $\pi_1$, \ie $\{p_i^{[t]}\}_{i}$, such that the polynomials $p_i^{[t]}$ have total degree less than or equal to $P_d$, is going to be denoted by
	\begin{equation}
	\label{eq:d2tot}
	\alpha_P^{[t]}(\omega)=\sum_{\|\pi_1^{-1}(i)\|_1\leq P_d,\ i\in\N} \tilde{c}_i^{[t]}p_i^{[t]}(\omega)=\sum_{i=0}^P \tilde{c}_i^{[t]}p_i^{[t]}(\omega),
	\end{equation}
	where $P+1={P_d+2 \choose 2}=\frac{(P_d+1)^2}{2}+\frac{(P_d+1)}{2}$. From \eqref{eq:piCantor}, we have $P=\pi_1(P_d,0)$.
		\item[Maximal Degree] is determined by the $\infty$-norm,  such that $\|(i_1,i_2)\|_\infty=\max \{i_1,i_2\}$. A pairing function, associated to the maximal degree, assigns consecutive numbers to points along the edges of squares of $\N\times\N$, \eg the Rosenberg-Strong pairing function $\pi_\infty(i_1,i_2)$ illustrated in the right pane of Figure~\ref{fig:pi} and specified by
	\begin{equation}\label{eq:piRosenberg}
	\pi_\infty(i_1,i_2)=(\max\{i_1,i_2\})^2+\max\{i_1,i_2\}+i_1-i_2.
	\end{equation}
	The polynomial approximation obtained from the corresponding basis $\{p_i^{[m]}\}_{i}$, such that the polynomials $p_i^{[m]}$ have maximal degree less than or equal to $P_d$ is going to be denoted by
	\begin{equation}
	\label{eq:d2max}
	\alpha_P^{[m]}(\omega)=\sum_{\|\pi_\infty^{-1}(i)\|_\infty\leq P_d,\ i\in\N} \tilde{c}_i^{[m]}p_i^{[m]}(\omega)=\sum_{i=0}^P \tilde{c}_i^{[m]}p_i^{[m]}(\omega),
	\end{equation}
	where $P+1=(P_d+1)^2$. From \eqref{eq:piRosenberg}, we have $P=\pi_\infty(P_d,0)$.
\end{description}

\subsection{Galerkin approach}\label{sec:gald2}
In this section, the spectral abscissa functions associated to 
\eqref{eq:spring} with control parameters of Table~\ref{tab:d2}, are approximated by the Galerkin approach  on total and maximal Legendre tensor basis, $\{p_i^{[t]}\}_i$ and $\{p_i^{[m]}\}_i$  with $\rho(\omega)=1/4$ for $\omega\in[-1,1]^2$. Along the same line of Section~\ref{sec:Galerkin}, first of all we consider the approximation error of the Galerkin approach, where the coefficients are evaluated with small  numerical error (w.r.t.~the approximation error); then we approximate the coefficients, $\tilde{c}_i^M$, via integration methods and we analyze the numerical error. At the end of this section, some advice to set the different parameters of this approach are given, in such a way that the numerical error does not dominate the approximation error. 

\subsubsection{Approximation error}\label{sec:gapproxd2}

\begin{figure}[!h]
   	\begin{center}
   		\begin{subfigure}[h]{0.45\linewidth}
   			\centering
   			\begin{tikzpicture}

\begin{axis}[%
width=0.8\textwidth,
height=0.5\textwidth,
at={(0\textwidth,0\textwidth)},
scale only axis,
xmode=log,
xmin=1,
xmax=100,
xtick={  1,  10, 100},
xticklabels={1,10,100},
xminorticks=false,
xlabel={$P_d+1$},
xmajorgrids,
ymode=log,
ymin=1e-11,
ymax=5,
ytick={1e-11,1e-10,1e-09,1e-08,1e-07,1e-06,1e-05,1e-04,1e-03,1e-02,1e-01,1},
ylabel={$\|\alpha-\alpha_P^{[t]}\|_\infty/\|\alpha\|_\infty$},
yminorticks=false,
ymajorgrids,
axis background/.style={fill=white},
ylabel style = {font=\small},
xlabel style = {font=\small},
title style={font=\bfseries, yshift=-0.2cm},
title={\small{Total Degree}}
]
\addplot [color=mycolor1,solid,line width=2.0pt,forget plot]
  table[row sep=crcr]{%
1	0.699594130506051\\
2	0.0473827589125243\\
3	0.000629941465585764\\
4	1.73675398410227e-05\\
5	7.49015726893161e-07\\
6	3.67330883072919e-08\\
7	7.85282030359704e-10\\
8	3.58102340308662e-11\\
9	2.52181588440503e-12\\
10	1.25796484263752e-12\\
11	1.27001753609417e-12\\
12	1.2854337719108e-12\\
13	1.27310078325749e-12\\
14	1.26048749940752e-12\\
15	1.24451067319755e-12\\
16	1.21395849676094e-12\\
17	1.2035875744843e-12\\
18	1.19854226094431e-12\\
19	1.21535997274427e-12\\
20	1.21031465920428e-12\\
21	1.2209658766776e-12\\
22	1.20667082164762e-12\\
23	1.22685207580758e-12\\
24	1.20975406881095e-12\\
25	1.21648115353094e-12\\
26	1.21535997274427e-12\\
27	1.20975406881095e-12\\
28	1.20835259282762e-12\\
29	1.19770137535431e-12\\
30	1.21648115353094e-12\\
31	1.20975406881095e-12\\
32	1.21227672558095e-12\\
33	1.22769296139758e-12\\
34	1.24226831162422e-12\\
35	1.22713237100425e-12\\
36	1.22236735266093e-12\\
37	1.20526934566429e-12\\
38	1.2035875744843e-12\\
39	1.26020720421085e-12\\
40	1.34513664880068e-12\\
41	1.23161709415091e-12\\
42	1.25628307145753e-12\\
43	1.2862746575008e-12\\
44	1.32355391865739e-12\\
45	1.33868985927736e-12\\
46	1.38129472917061e-12\\
47	1.33014085577905e-12\\
48	2.18181781084904e-12\\
49	1.91918121157289e-12\\
50	1.61029590484683e-12\\
51	1.78912424031981e-12\\
52	2.53134592109168e-12\\
53	2.2151729392523e-12\\
54	1.4070818872639e-12\\
55	1.42053605670387e-12\\
56	1.62767420704013e-12\\
57	2.1668220178274e-12\\
58	2.03101899504267e-12\\
59	2.34999492884871e-12\\
60	1.49705664539372e-12\\
61	1.62823479743346e-12\\
62	1.63159833979345e-12\\
63	1.74483759924656e-12\\
64	1.71848985075995e-12\\
65	1.69186180707667e-12\\
66	1.87069014254965e-12\\
67	2.15883360472242e-12\\
68	3.18751697648706e-12\\
69	4.08922662416195e-12\\
70	3.14967712493713e-12\\
71	2.00523183694938e-12\\
72	3.04919129693233e-12\\
73	1.71470586560496e-12\\
74	2.00943626489938e-12\\
75	1.90096202378959e-12\\
76	1.93992305612618e-12\\
77	2.66504672990142e-12\\
78	2.90792251781261e-12\\
79	3.80851098470083e-12\\
80	3.83163533842579e-12\\
81	2.17396954534239e-12\\
82	1.98000526924943e-12\\
83	2.97981823575747e-12\\
84	3.32934634600011e-12\\
85	3.217788857727e-12\\
86	3.59394501165292e-12\\
87	2.18826460037236e-12\\
88	2.17144688857239e-12\\
89	1.9786037932661e-12\\
90	4.12790736130187e-12\\
91	4.33056078849147e-12\\
92	2.72138606443131e-12\\
93	3.22591741843031e-12\\
94	2.92249786803925e-12\\
95	2.64850931329812e-12\\
96	2.42819728871855e-12\\
97	2.90694148462428e-12\\
98	2.39372097952862e-12\\
99	2.09366497149754e-12\\
100	3.16369188477044e-12\\
101	2.78361159809119e-12\\
102	2.39063773236529e-12\\
};
\addplot [color=mycolor2,solid,line width=2.0pt,forget plot]
  table[row sep=crcr]{%
1	0.824878104407754\\
2	0.875111890725568\\
3	0.564069022105175\\
4	0.683845116463415\\
5	0.43768997226538\\
6	0.407621905385209\\
7	0.402162461094579\\
8	0.30488827503339\\
9	0.394177541809577\\
10	0.277059129516743\\
11	0.298843402409798\\
12	0.192932993609224\\
13	0.165179057192495\\
14	0.240091726365892\\
15	0.183096042924442\\
16	0.223941021463726\\
17	0.156677234452898\\
18	0.1623017689597\\
19	0.143913141554699\\
20	0.121346623972162\\
21	0.165986472049535\\
22	0.125517372239288\\
23	0.147390155074719\\
24	0.100770956495253\\
25	0.0994412095858696\\
26	0.116024042021078\\
27	0.0955470234541898\\
28	0.124803171490983\\
29	0.0923929510074877\\
30	0.10448539750283\\
31	0.0724794176529657\\
32	0.0671147287081508\\
33	0.0978880684316062\\
34	0.078327007266745\\
35	0.098061315974568\\
36	0.0702581808594133\\
37	0.0759614841623577\\
38	0.0690739947971854\\
39	0.0609312174717252\\
40	0.0844794591120475\\
41	0.0655354185640553\\
42	0.0790287907698375\\
43	0.0545229598262065\\
44	0.0558010532715803\\
45	0.0670451743528179\\
46	0.0566439352967526\\
47	0.0750194677287156\\
48	0.0560602079071369\\
49	0.0646297910975394\\
50	0.0472337551070329\\
51	0.0444589004517891\\
52	0.0649371527719373\\
53	0.0522350057318011\\
54	0.0657317212529289\\
55	0.046926448320023\\
56	0.0508446953225267\\
57	0.0501298777328578\\
58	0.0443029374429357\\
59	0.0609393957252751\\
60	0.0472692138922115\\
61	0.0569447108967115\\
62	0.0388493355044502\\
63	0.0395416320882873\\
64	0.0519180388969136\\
65	0.0436395383763018\\
66	0.0573692716228263\\
67	0.0425003814876974\\
68	0.0486639242799197\\
69	0.0391551062517988\\
70	0.0366924163991903\\
71	0.0527706689471522\\
72	0.0423452629150985\\
73	0.0531030616458677\\
74	0.0373772392514847\\
75	0.0401549018111079\\
76	0.043854639286337\\
77	0.0383671753397516\\
78	0.0523219562740875\\
79	0.0399651455239176\\
80	0.0477371229387781\\
81	0.0323974921260145\\
82	0.0319576175806817\\
83	0.0477626849552527\\
84	0.0394882133220129\\
85	0.0511838249852012\\
86	0.0370437321399409\\
87	0.041637781358055\\
88	0.0398682399609384\\
89	0.036168635799101\\
90	0.0508093722798476\\
91	0.0397558627832681\\
92	0.0489130282597984\\
93	0.0330692498135114\\
94	0.034172277968835\\
95	0.0467354913988231\\
96	0.0395692578221786\\
97	0.052496033535784\\
98	0.0386873150376938\\
99	0.0447103986303231\\
100	0.0391345433493112\\
101	0.0366009771950212\\
102	0.0525280685326862\\
};
\addplot [color=mycolor3,solid,line width=2.0pt,forget plot]
  table[row sep=crcr]{%
1	2.17797143814889\\
2	2.01891941426668\\
3	1.60470793584719\\
4	1.60359873125033\\
5	1.40971574457722\\
6	1.40834849173188\\
7	1.28041678741419\\
8	1.27805377208898\\
9	1.18560288171171\\
10	1.18319834473134\\
11	1.11206328194906\\
12	1.11033419340514\\
13	1.05248669745584\\
14	1.05126841834758\\
15	1.00365410496766\\
16	1.00250269507742\\
17	0.963276392763444\\
18	0.962078625273276\\
19	0.928833231027028\\
20	0.92765928116374\\
21	0.898298138859444\\
22	0.897201088348358\\
23	0.870969089440347\\
24	0.869989580921912\\
25	0.846778934112667\\
26	0.845958488624858\\
27	0.825394876378324\\
28	0.824774635574548\\
29	0.806361815464627\\
30	0.805911890333825\\
31	0.789103784021212\\
32	0.789091958408642\\
33	0.773566521330775\\
34	0.774277706766201\\
35	0.759554190504849\\
36	0.761076902278528\\
37	0.746541452272905\\
38	0.748681424020278\\
39	0.734064897224471\\
40	0.736329226712062\\
41	0.721612446828567\\
42	0.723810861128127\\
43	0.709169867801932\\
44	0.711451021705875\\
45	0.696999717460955\\
46	0.699416215200804\\
47	0.685120510502453\\
48	0.687679301736059\\
49	0.673744893648927\\
50	0.676527942342213\\
51	0.663253803284921\\
52	0.666446565100501\\
53	0.660616611747246\\
54	0.663755498821268\\
55	0.658400987507739\\
56	0.661364440531254\\
57	0.655542700980717\\
58	0.658190251305244\\
59	0.652239885180819\\
60	0.654628575520248\\
61	0.648949129780947\\
62	0.651163552908769\\
63	0.647280404307093\\
64	0.649951370110671\\
65	0.652290058917263\\
66	0.656223236285281\\
67	0.656806089783043\\
68	0.660091342280706\\
69	0.660293627031405\\
70	0.662833096006075\\
71	0.665141939928548\\
72	0.669284419756548\\
73	0.672544075216059\\
74	0.679476665789054\\
75	0.681104121899188\\
76	0.687185343839508\\
77	0.688983068768424\\
78	0.693818591974294\\
79	0.697766220064493\\
80	0.704858122343288\\
81	0.708081712385113\\
82	0.715882168604149\\
83	0.717920694211\\
84	0.725260844428869\\
85	0.728274977077557\\
86	0.736613314423476\\
87	0.738958326670665\\
88	0.747264458764259\\
89	0.749525337574786\\
90	0.75780945307552\\
91	0.761008443269228\\
92	0.77127305247002\\
93	0.772328584925309\\
94	0.781150764904678\\
95	0.78297114834215\\
96	0.792052018312552\\
97	0.793389230008802\\
98	0.801952924856296\\
99	0.803550992452249\\
100	0.812018205806837\\
101	0.813977307284238\\
102	0.823333427497395\\
};
\end{axis}
\end{tikzpicture}%
   		\end{subfigure}\quad
   		\begin{subfigure}[h]{0.45\linewidth}
   			\centering
   			\begin{tikzpicture}

\begin{axis}[%
width=0.8\textwidth,
height=0.5\textwidth,
at={(0\textwidth,0\textwidth)},
scale only axis,
xmode=log,
xmin=1,
xmax=100,
xlabel={$P_d+1$},
xtick={  1,  10, 100},
xticklabels={1,10,100},
xminorticks=false,
xmajorgrids,
ymode=log,
ymin=1e-11,
ymax=5,
yminorticks=false,
ytick={1e-11,1e-10,1e-09,1e-08,1e-07,1e-06,1e-05,1e-04,1e-03,1e-02,1e-01,1},
ylabel={$\|\alpha-\alpha_P^{[m]}\|_\infty/\|\alpha\|_\infty$},
ymajorgrids,
axis background/.style={fill=white},
ylabel style = {font=\small},
xlabel style = {font=\small},
title style={font=\bfseries, yshift=-0.2cm},
title={\small{Maximal Degree}}
]
\addplot [color=mycolor1,solid,line width=2.0pt,forget plot]
  table[row sep=crcr]{%
1	0.699594130506046\\
2	0.00348032929494438\\
3	3.81829572075556e-05\\
4	3.53486573831954e-07\\
5	1.62862435169779e-08\\
6	1.71211874093973e-10\\
7	4.37933215271138e-12\\
8	1.30841797803742e-12\\
9	1.29552439899078e-12\\
10	1.29692587497411e-12\\
11	1.28991849505746e-12\\
12	1.28739583828747e-12\\
13	1.25376041468753e-12\\
14	1.24535155878755e-12\\
15	1.23077620856091e-12\\
16	1.23105650375758e-12\\
17	1.24899539634421e-12\\
18	1.24899539634421e-12\\
19	1.23694270288756e-12\\
20	1.24451067319755e-12\\
21	1.24198801642755e-12\\
22	1.25544218586753e-12\\
23	1.24675303477088e-12\\
24	1.2506771675242e-12\\
25	1.24030624524756e-12\\
26	1.24198801642755e-12\\
27	1.2349806365109e-12\\
28	1.23161709415091e-12\\
29	1.23638211249423e-12\\
30	1.24815451075421e-12\\
31	1.23554122690423e-12\\
32	1.24675303477088e-12\\
33	1.24675303477088e-12\\
34	1.2515180531142e-12\\
35	1.27590373522415e-12\\
36	1.29117982344246e-12\\
37	1.34934107675068e-12\\
38	1.35915140863399e-12\\
39	1.29286159462245e-12\\
40	1.28935790466413e-12\\
41	1.27786580160082e-12\\
42	1.28907760946746e-12\\
43	1.33168247936071e-12\\
44	1.42081635190054e-12\\
45	1.40315775451057e-12\\
46	1.36307554138732e-12\\
47	1.29832735095744e-12\\
48	1.31066033961075e-12\\
49	1.38858240428393e-12\\
50	1.65121900356008e-12\\
51	1.42698284622719e-12\\
52	1.38045384358062e-12\\
53	1.45473207069714e-12\\
54	1.45921679384379e-12\\
55	1.51555612837368e-12\\
56	1.4202557615072e-12\\
57	1.47603450564376e-12\\
58	1.42221782788387e-12\\
59	1.49929900696705e-12\\
60	1.56825162534691e-12\\
61	1.53769944891031e-12\\
62	1.48332218075708e-12\\
63	1.47351184887377e-12\\
64	1.58198608998355e-12\\
65	1.81463110321643e-12\\
66	1.4227784182772e-12\\
67	1.53797974410697e-12\\
68	1.40792277285389e-12\\
69	1.58562992754021e-12\\
70	1.43006609339052e-12\\
71	1.46902712572711e-12\\
72	1.50112092574538e-12\\
73	1.49845812137705e-12\\
74	1.52004085152034e-12\\
75	1.84350150847304e-12\\
76	1.90250364737125e-12\\
77	1.65514313631341e-12\\
78	1.72185339311994e-12\\
79	1.73110313460992e-12\\
80	1.94188512250284e-12\\
81	1.69074062629e-12\\
82	1.87265220892631e-12\\
83	2.01616334961936e-12\\
84	1.79949516259646e-12\\
85	1.55830114586526e-12\\
86	1.71512630839996e-12\\
87	1.75801147348987e-12\\
88	2.12870187108081e-12\\
89	2.41614459526191e-12\\
90	3.40166250673997e-12\\
91	3.78524648337755e-12\\
92	2.56820473945328e-12\\
93	2.62622584516316e-12\\
94	3.99266492891047e-12\\
95	3.03854007945902e-12\\
96	1.80412003334145e-12\\
97	2.02148895835602e-12\\
98	1.89956054780626e-12\\
99	2.24656600127891e-12\\
100	2.21881677680896e-12\\
101	2.31187478210211e-12\\
102	2.60226060584821e-12\\
};
\addplot [color=mycolor2,solid,line width=2.0pt,forget plot]
  table[row sep=crcr]{%
1	0.824877840638045\\
2	0.559991742941474\\
3	0.42247221810307\\
4	0.334322692528679\\
5	0.28292381397827\\
6	0.245677173285718\\
7	0.217834574156964\\
8	0.192902959315743\\
9	0.166332857243716\\
10	0.139381267098198\\
11	0.11617938451153\\
12	0.105898923373164\\
13	0.0980998711348993\\
14	0.0906275238148083\\
15	0.0838913483859985\\
16	0.0776082124728613\\
17	0.0729378480220278\\
18	0.0685430075836753\\
19	0.0647540603481319\\
20	0.0612042860679443\\
21	0.0578521756340246\\
22	0.0550615537754452\\
23	0.0524118244876906\\
24	0.0502573165141479\\
25	0.048085919611791\\
26	0.0459981273458753\\
27	0.0440599857562174\\
28	0.0425521832693805\\
29	0.0410174864052563\\
30	0.0396040290761785\\
31	0.0380800529726638\\
32	0.0367201665605631\\
33	0.035813329885715\\
34	0.0346063813802255\\
35	0.0336787907035905\\
36	0.0324328488845833\\
37	0.0316664243129843\\
38	0.0308377849193902\\
39	0.0299027057928713\\
40	0.0292920540719827\\
41	0.028476324410997\\
42	0.0278141302471579\\
43	0.0270806401918475\\
44	0.0262129729855494\\
45	0.026028232869079\\
46	0.0254335066828857\\
47	0.0248574109275374\\
48	0.0239982397729764\\
49	0.0236847872227549\\
50	0.0236119913558\\
51	0.0231391990232826\\
52	0.0223202777123465\\
53	0.0219703226166461\\
54	0.0218353877975474\\
55	0.0216799843056336\\
56	0.0212783381678664\\
57	0.0206339600898737\\
58	0.0204931832444606\\
59	0.020429479869253\\
60	0.0200719068038725\\
61	0.0197988060002329\\
62	0.0193538878097721\\
63	0.0193965410710114\\
64	0.0191344938609425\\
65	0.0187983038990901\\
66	0.0186001577881644\\
67	0.0185255895657694\\
68	0.0183311025911728\\
69	0.0178915870475755\\
70	0.01777206302715\\
71	0.0177436795752898\\
72	0.0176190156503911\\
73	0.017310830391504\\
74	0.0170825554756655\\
75	0.0170448079385496\\
76	0.0170328606247793\\
77	0.016828623936842\\
78	0.0165293975784858\\
79	0.016482489347738\\
80	0.0165618374253173\\
81	0.016354904425166\\
82	0.0160429959459485\\
83	0.0160288022699725\\
84	0.0160912463262697\\
85	0.0159271195539024\\
86	0.0156749892835452\\
87	0.0156418762054733\\
88	0.015649135371508\\
89	0.0155919870723935\\
90	0.0153020147698479\\
91	0.0149878976306419\\
92	0.0149714660609149\\
93	0.0148538356027369\\
94	0.0146100402318169\\
95	0.0144138567894736\\
96	0.0144970810008086\\
97	0.0146516157137735\\
98	0.0145494904120972\\
99	0.0145445924460283\\
100	0.0146341851426189\\
101	0.0146388108172621\\
102	0.0147618589970755\\
};
\addplot [color=mycolor3,solid,line width=2.0pt,forget plot]
  table[row sep=crcr]{%
1	2.17800331823441\\
2	2.01863828738248\\
3	1.47221378594913\\
4	1.4594610109565\\
5	1.24881762213803\\
6	1.23743369181826\\
7	1.11669917104835\\
8	1.10635159970639\\
9	1.02542181203791\\
10	1.01588497615614\\
11	0.956958789863775\\
12	0.948003763208917\\
13	0.902739227926588\\
14	0.894228000897145\\
15	0.858216973483798\\
16	0.850065029804217\\
17	0.820805487522576\\
18	0.812873626548509\\
19	0.788556173885151\\
20	0.780805700239636\\
21	0.76027551365529\\
22	0.752677275163105\\
23	0.735176921949892\\
24	0.727702231391316\\
25	0.712590153055274\\
26	0.705244070160141\\
27	0.692140811567058\\
28	0.684899947155668\\
29	0.673439988673538\\
30	0.666357214288613\\
31	0.65632425463006\\
32	0.64936638968722\\
33	0.64049453476354\\
34	0.633753822009652\\
35	0.625980294749576\\
36	0.619406776530315\\
37	0.612583131873677\\
38	0.606211973091595\\
39	0.600060725133365\\
40	0.593965766016731\\
41	0.58846006275038\\
42	0.582486887934043\\
43	0.57738709545035\\
44	0.571435103810384\\
45	0.566809470408435\\
46	0.560787984305744\\
47	0.556561215386419\\
48	0.550583792054235\\
49	0.546759790371408\\
50	0.540758740889297\\
51	0.537170726399008\\
52	0.531211999300948\\
53	0.52790036968505\\
54	0.522007728487229\\
55	0.51897340198599\\
56	0.512930196177884\\
57	0.509994059157672\\
58	0.504258403365196\\
59	0.501673868161825\\
60	0.495919579798562\\
61	0.493557666434463\\
62	0.487947076515542\\
63	0.485858496053465\\
64	0.480532866883953\\
65	0.478863200769902\\
66	0.473375519203453\\
67	0.471923430348465\\
68	0.466320199620021\\
69	0.465060046836309\\
70	0.459715620628885\\
71	0.458762433623928\\
72	0.453343652589912\\
73	0.452638871924385\\
74	0.449401518751666\\
75	0.447471657577058\\
76	0.447342644859616\\
77	0.445007527494016\\
78	0.445157586657196\\
79	0.442204378217866\\
80	0.442571768857644\\
81	0.438990250981174\\
82	0.439592141630804\\
83	0.435510482081875\\
84	0.436493872740277\\
85	0.432146505192383\\
86	0.433600701493204\\
87	0.42928579999462\\
88	0.431317618483184\\
89	0.426975159700467\\
90	0.429409177940693\\
91	0.424791244820438\\
92	0.427301836628261\\
93	0.422536323967625\\
94	0.424857663259717\\
95	0.420087661978392\\
96	0.422130234054379\\
97	0.41882549949216\\
98	0.42096887403604\\
99	0.421226411358972\\
100	0.422856099959876\\
101	0.423724556221337\\
102	0.424918114428832\\
};
\end{axis}
\end{tikzpicture}%
   		\end{subfigure}
   		\vspace{-0.01\textwidth}  
   		\begin{tikzpicture}
   		\hspace{0.05\linewidth}
   		\begin{customlegend}[legend columns=5,legend style={align=left,draw=none,column sep=1ex},legend entries={\sae{SAE},\mssae{MSSAEs},\mnssae{MNSSAEs}}]
   		\addlegendimage{color=mycolor1, mark=none,solid,line width=2.0pt,line legend}
   		\addlegendimage{color=mycolor2, mark=none,solid,line width=2.0pt}   
   		\addlegendimage{color=mycolor3, mark=none,solid,line width=2.0pt}
   		\end{customlegend}
   		\end{tikzpicture}
   	\end{center}
   	\caption{Approximation errors of the Galerkin approach to compute the polynomial approximation w.r.t.~total degree $\alpha_P^{[t]}$ and maximal degree $\alpha_P^{[m]}$,  for the benchmark in Example~\ref{ex:d2}.} 
   	\label{fig:gd2_approx_err}
   \end{figure}
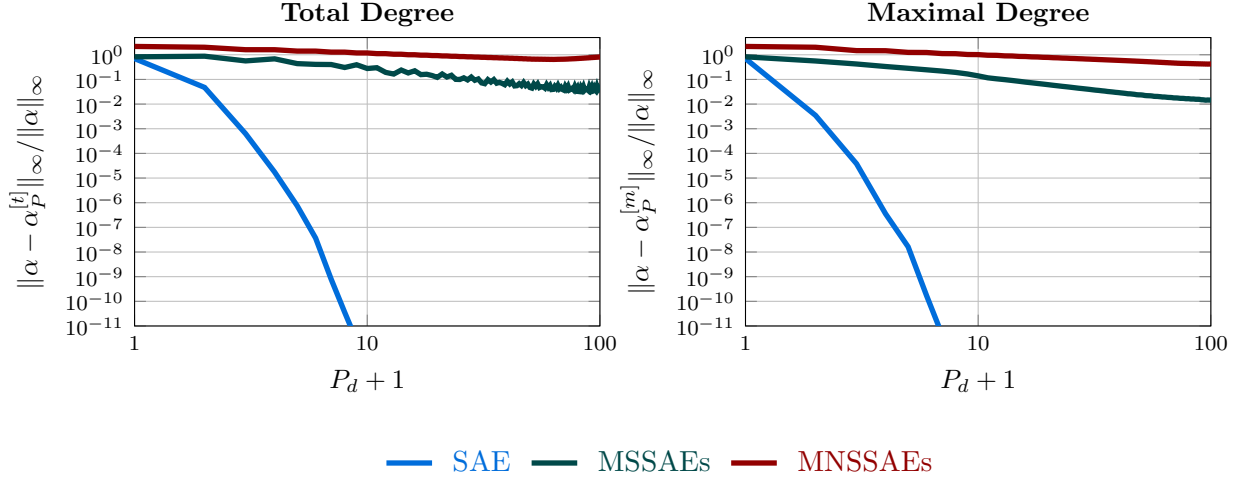

 Figure~\ref{fig:gd2_approx_err} shows the convergence rates of the approximation error with the Galerkin approach, where the coefficients are evaluated with a numerical error which does not affect the approximation error, following the advice of the upcoming Section~\ref{sec:coupling}. For both total and maximal degrees, the converge rates are analogous, and they are approximately of order $\mathcal{O}(P_d^{-P_d})$,  $\mathcal{O}(P_d^{-1})$ and   $\mathcal{O}(P_d^{-0.3})$, for the  \sae{SAE}, \mssae{MSSAEs}, and \mnssae{MNSSAEs} cases in Example~\ref{ex:d2}, respectively.%
   
\subsubsection{Numerical error}\label{sec:gnumd2}
The integral computation, needed to evaluate the coefficients in the Galerkin approach (left hand side of \eqref{eq:ci}), is done, in this section, by bi-dimensional generalizations of the Clenshaw-Curtis quadrature rule based on Chebyshev points. The integration methods, here considered, are based on the following set of points:
\begin{description}
	\item[Tensor product Chebyshev grid] is constructed by tensor product of $(m_c+1)$ unidimensional Chebyshev points. The total number of points used is $M+1=(m_c+1)^2$.
	\item[Padua points] are the self-intersections and boundary contacts of the generating curve $T_{m_p}(x)-T_{m_p+1}(y)=0$ for $(x,y)\in[-1,1]^2$, where $T_{m_p}$ is the Chebyshev polynomial of degree $m_p$  (cf., \eg, \cite{Bos2006, Caliari2008} and references therein). The total number of points is $M+1={m_p+2\choose 2}$.
\end{description}

The approximation of an integral, obtained by integrating the interpolant of the integrand, evaluated on the previous sets of points, leads to the following integration methods (cf. \eg \cite{Sommariva2008})
\begin{description}
	\item[Tensorial Clenshaw–-Curtis cubature rule] based on tensor product Chebyshev grid.
	\item[Non-tensorial Clenshaw–-Curtis cubature rule] relying on Padua points. 
\end{description} 
These methods approximate the coefficients $\tilde{c}_i^{[t],M}$ and $\tilde{c}_i^{[m],M}$ and, hence, permit to compute the polynomial approximations w.r.t.~total and maximal degrees, \ie $\alpha_P^{[t],M}$ and $\alpha_P^{[m],M}$, respectively. By the advice furnished in the upcoming Section~\ref{sec:coupling}, the coefficients $\tilde{c}_i^{[t],M}$ are computed by non-tensorial Clenshaw-Curtis cubature rules on $M+1$ Padua points, while $\tilde{c}_i^{[m],M}$ are evaluated by tensorial Clenshaw-Curtis cubature rules on $M+1$ tensor product Chebyshev grid. 

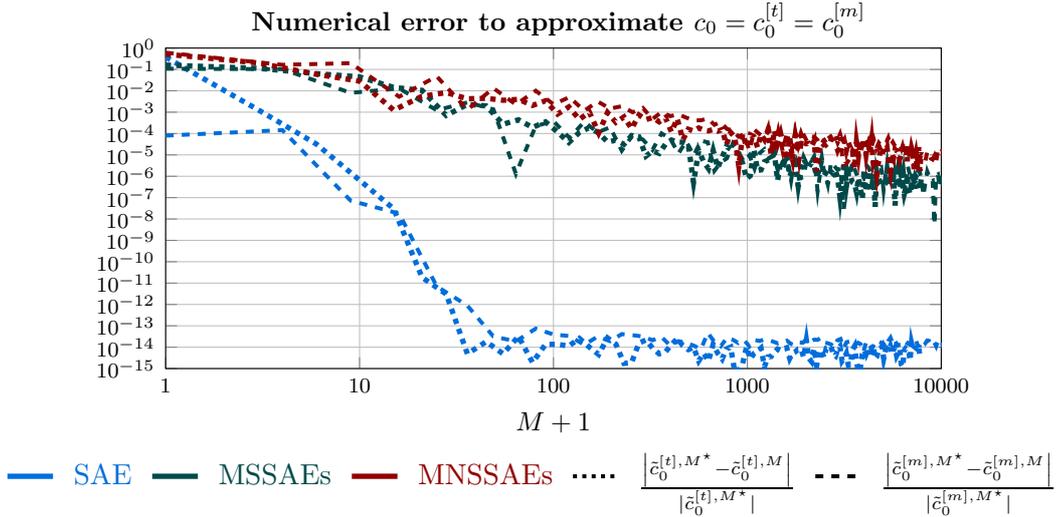
\begin{figure}[!h]
	\begin{center}
		\include{c0d2}
		\vspace{-0.03\textwidth}  
		\begin{tikzpicture}
		\begin{customlegend}[legend columns=5,legend style={align=left,draw=none,column sep=1ex},legend entries={\sae{SAE},\mssae{MSSAEs},\mnssae{MNSSAEs},\small{$\frac{\left\lvert\tilde{c}_0^{[t],M^\star}-\tilde{c}_0^{[t],M}\right\rvert}{\lvert\tilde{c}_0^{[t],M^\star}\rvert}$}, \small{$\frac{\left\lvert\tilde{c}_0^{[m],M^\star}-\tilde{c}_0^{[m],M}\right\rvert}{\lvert\tilde{c}_0^{[m],M^\star}\rvert}$} }]
		\addlegendimage{color=mycolor1, mark=none,solid,line width=2.0pt,line legend}
		\addlegendimage{color=mycolor2, mark=none,solid,line width=2.0pt}   
		\addlegendimage{color=mycolor3, mark=none,solid,line width=2.0pt}
		\addlegendimage{color=black, mark=none,dotted,line width=1.5pt}
		\addlegendimage{color=black, mark=none,dashed,line width=1.5pt}
		\end{customlegend}
		\end{tikzpicture}
		\caption{Numerical error to compute the first coefficient ${c}_0$ of \eqref{eq:spapproxd2} for benchmark in Example~\ref{ex:d2}. The approximations $\tilde{c}_0^{[t],M}$, and $\tilde{c}_0^{[m],M}$ are obtained by non-tensorial and tensorial Clenshaw-Curtis cubature rules, respectively. The reference values are computed by $M^\star>5\cdot10^{5}$ which corresponds to $m_p=999$ and $m_c=707$.} 
		\label{fig:c0d2}
	\end{center}
\end{figure}

In order to compute the numerical errors, $\tilde{c}_i^{[t],M^\star}$ and $\tilde{c}_i^{[m],M^\star}$ are considered as reference values. These coefficients are computed by the corresponding cubature rules based on $M^\star>5\cdot10^{5}$  points (in particular $m_p=999$ and $m_c=707$). Figure~\ref{fig:c0d2} shows the error to compute the first coefficient $\tilde{c}_0^M$, which is independent from the multivariate polynomial degree, since ${c}_0^{[t]}={c}_0^{[m]}$.  The cubature rules on tensor product Chebyshev grid and Padua points present similar convergences rates; in particular  the \mssae{MSSAEs} and \mnssae{MNSSAEs} converge almost as $\mathcal{O}(M^{-2})$ and $\mathcal{O}(M^{-1.5})$, respectively; while the \sae{SAE} present a convergence of order $\mathcal{O}(M^{-M})$.

In terms of a PC expansion (see Theorem~\ref{th:pce}), Figure~\ref{fig:c0d2}  corresponds to the error of the mean of the spectral abscissa $\alpha(\bomega)$ associated to 
\eqref{eq:spring}, where $\omega$ is a realization of $\bomega$, random vector uniformly distributed in $\Su$. For this example, the tensorial and non-tensorial Clenshaw-Curtis cubature rules provide a faster convergence rate  of the mean w.r.t.~Quasi Monte Carlo method based on Halton sequences, which converges with an order of $\mathcal{O}(M^{-1})$ for all the three cases (\sae{SAE}, and \mssae{MSSAEs} and \mnssae{MNSSAEs})  as shown in the left pane of Figure 4 in \cite{Fenzi2017}.

\subsubsection{Decoupling numerical error w.r.t.~approximation error}\label{sec:coupling}
This section furnishes some advice on how to set the number of points for cubature rules. We should consider integration method based on a small number of points, $M$, which are, however, still sufficiently large in order that the approximation error in the Galerkin approach is not corrupted by the numerical error due to the approximation of the integrals. 

Since the coefficients $c_i$ in \eqref{eq:spabsapprox} and \eqref{eq:spapproxd2} are determined by integrating $\alpha(\omega)p_i(\omega)\rho(\omega)$ over the domain $\omega\in\Su$, we request that the polynomial approximation behind the integration rule is exact for all polynomials $\{p_i\}_{i=0}^P$. 

This advice leads in the unidimensional case (Section~\ref{sec:gd1num}) that the number of nodes of the interpolatory quadrature rules (\ie Gauss and Clenshaw-Curtis quadrature rules), $M+1$, is greater than the degree of the polynomial approximation $\alpha_P^M$, \ie $P$.

For the bi-dimensional case, the polynomial approximation behind the non-tensorial Clenshaw-Curtis cubature rule, \ie the interpolation on Padua points, is exact for  polynomial with total degree less than or equal to $m_p$. The interpolant on tensor product Chebyshev grid, underling the tensorial  Clenshaw-Curtis cubature rule, is exact for all the polynomials with maximal degree less than or equal to $m_c$. From this properties we can observe that the interpolatory cubature rule for $D>1$ are associated to the polynomial multivariate degree. In particular, non-tensorial and tensorial Clenshaw-Curtis cubature rules are interpolatory cubature rules associated to the total and maximal degrees for $D=2$, respectively.

 Following the advice, $\alpha_P^{[t],M}$ is accurately approximated by non-tensorial Clenshaw-Curtis cubature rule with $m_p\geq P_d$ or by tensorial Clenshaw-Curtis cubature rule with $m_c\geq P_d$. On the other hand, $\alpha_P^{[m],M}$ can be evaluated   by tensorial Clenshaw-Curtis with $m_c\geq P_d$  or by non-tensorial Clenshaw-Curtis cubature rule with $m_p\geq 2 P_d$. Therefore, the number of points used in the integration rule are minimized if we consider non-tensorial Clenshaw-Curtis cubature rule  for approximating the coefficients of $\alpha_P^{[t],M}$  with $m_p\geq P_d$, and tensorial  Clenshaw-Curtis cubature rule for approximating the coefficients of $\alpha_P^{[m],M}$  with $m_c\geq P_d$.

 \begin{table}[!h]
 	\caption{Advice to compute Galerkin polynomial approximations \eqref{eq:d2tot} and \eqref{eq:d2max}, whose coefficients are evaluated by non-tensorial and tensorial Clenshaw-Curtis cubature rule, based on Padua points and tensor product Chebyshev grid, respectively. The natural choices are highlighted.}
 	\label{tab:decoupling}
 	\begin{center}
 		\begin{tabular}{lcc}
 			&\textbf{Total degree} $\alpha_P^{[t],M}$ & \bfseries{Maximal degree} $\alpha_P^{[m],M}$\\
 			\addlinespace[0.2em]
 			\multirow{ 2}{*}{\bfseries{Padua points}} 						& \cellcolor{celltab} $m_p\geq P_d$  	& $m_p\geq 2\cdot P_d$ \\ 
 			&\cellcolor{celltab} {\footnotesize$M+1={m_p+2 \choose 2}\geq P+1={P_d+2 \choose 2}$} & {\footnotesize $M+1={2\cdot m_p+2 \choose 2}\gg P+1=(P_d+1)^2$ }\\
 			\addlinespace[0.2em]
 			\bfseries{Tensor product}  	& $m_c\geq P_d$	&\cellcolor{celltab} $m_c\geq P_d$\\ 
 			\bfseries{Chebyshev grid}	& {\footnotesize$M+1=(m_c+1)^2\gg P+1={P_d+2 \choose 2}$ }&\cellcolor{celltab} {\footnotesize $M+1=(m_c+1)^2\geq P+1={(P_d+1)^2}$}\\
 		\end{tabular}
 	\end{center}
 \end{table}

  \begin{figure}[!h]
  	\begin{center}
  		\begin{subfigure}[h]{0.45\linewidth}
  			\centering
  			\include{decD2_total}
  			\vspace{-0.08\textwidth}  
  			\begin{tikzpicture}
  			\begin{customlegend}[legend columns=2,legend style={align=left,draw=none,column sep=1ex},legend entries={\small{$\sum_i\left\lvert c_i^{[t],M^\star}-\tilde{c}_i^{[t],M}\right\rvert$}, \small{$\frac{\|\alpha-\alpha_P^{[t],M}\|_\infty}{\|\alpha\|_\infty}$}}]
  			\addlegendimage{color=black, mark=none,dashdotted,line width=1.0pt}
  			\addlegendimage{color=black, mark=none,solid,line width=1.5pt}
  			\end{customlegend}
  			\end{tikzpicture}
  		\end{subfigure}\quad
  		\begin{subfigure}[h]{0.45\linewidth}
  			\centering
  			\include{decD2_max}
  			\vspace{-0.08\textwidth}  
  			\begin{tikzpicture}
  			\hspace{0.05\linewidth}
  			\begin{customlegend}[legend columns=2,legend style={align=left,draw=none,column sep=1ex},legend entries={\small{$\sum_i\left\lvert c_i^{[m],M^\star}-\tilde{c}_i^{[m],M}\right\rvert$}, \small{$\frac{\|\alpha-\alpha_P^{[m],M}\|_\infty}{\|\alpha\|_\infty}$}}]
  			\addlegendimage{color=black, mark=none,dashdotted,line width=1.0pt}
  			\addlegendimage{color=black, mark=none,solid,line width=1.5pt}
  			\end{customlegend}
  			\end{tikzpicture}
  		\end{subfigure}
  		\vspace{-0.0\textwidth}  
  		\begin{tikzpicture}
  		\begin{customlegend}[legend columns=3,legend style={align=left,draw=none,column sep=1ex},legend entries={\sae{SAE},\mssae{MSSAEs},\mnssae{MNSSAEs}}]
  		\addlegendimage{color=mycolor1, mark=none,solid,line width=2.0pt,line legend}
  		\addlegendimage{color=mycolor2, mark=none,solid,line width=2.0pt}   
  		\addlegendimage{color=mycolor3, mark=none,solid,line width=2.0pt}
  		\end{customlegend}
  		\end{tikzpicture}
  	\end{center}
  	\caption{Numerical and approximation errors to compute the polynomial approximation  by the Galerkin approach of the spectral abscissa functions of Example~\ref{ex:d2}, using the two non-natural choices of Table~\ref{tab:decoupling}.   The coefficients $\tilde{c}_i^{[t],M}$ and $\tilde{c}_i^{[m],M}$ are computed by tensorial and non-tensorial Clenshaw-Curtis cubature rules, respectively, where $M>P$ is given by $m_c=71<100$ and $m_p=140<2\cdot100$.} 
  	\label{fig:countercheck}
  \end{figure}
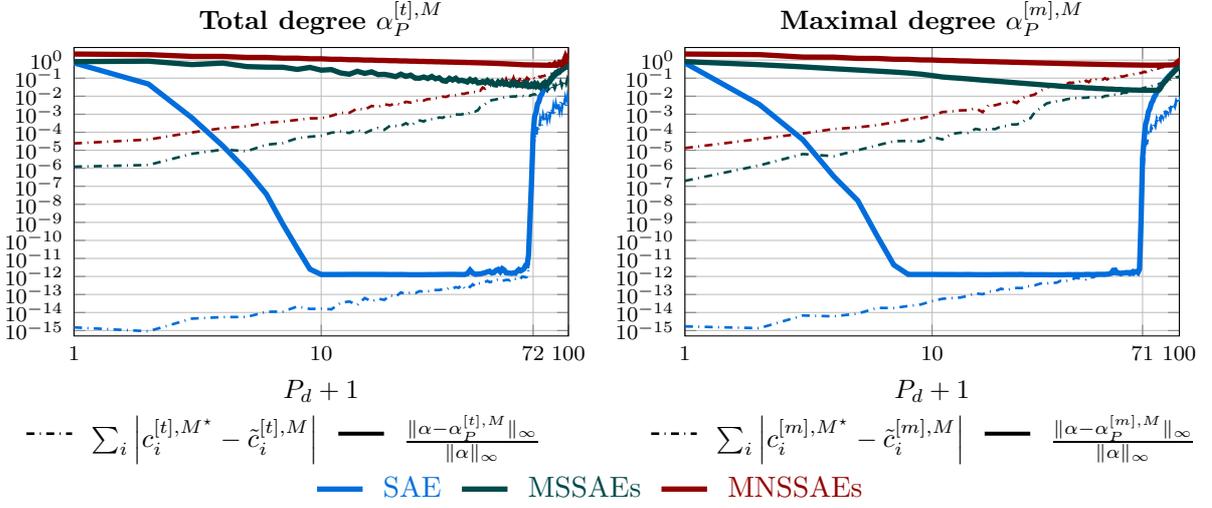

 Table~\ref{tab:decoupling} summarizes the advice for the bi-dimensional case, highlighting the natural choices of computing the polynomial approximations $\alpha_P^{[t],M}$ and $\alpha_P^{[m],M}$ through the associated interpolatory cubature rules. 
 These natural choices are used to obtain the convergence errors in Figure~\ref{fig:gd2_approx_err}, where the coefficients of $\alpha_P^{[t]}$,  $\alpha_P^{[m]}$ are evaluated with $m_p=110$ non-tensorial and $m_c=110$ tensorial Clenshaw-Curtis cubature rule, respectively. On the other hand, Figure~\ref{fig:countercheck} shows a counter check of this advice for spectral abscissa polynomial approximations of Example~\ref{ex:d2}. The number of points of the cubature rules $M$ are greater than the number of coefficients,  $P$,  that we want to approximate. However, the polynomial approximation $\alpha_P^{[t],M}$ is computed though tensorial Clenshaw-Curtis cubature with $m_c=71$ while the polynomial approximation $\alpha_P^{[m],M}$ is evaluated by non-tensorial Clenshaw-Curtis cubature rule with $m_p=140$. The numerical errors,
 \begin{equation}
 \sum_{\|\pi_1^{-1}(i)\|_1\leq P_d}\left\lvert c_i^{[t],M^\star}-\tilde{c}_i^{[t],M}\right\rvert\quad\text{and} \quad \sum_{\|\pi_\infty^{-1}(i)\|_\infty\leq P_d}\left\lvert c_i^{[m],M^\star}-\tilde{c}_i^{[m],M}\right\rvert,
  \end{equation} heavily affects the approximation error as soon as $P_d \gtrapprox m_c=71$ for $\alpha_P^{[t],M}$, and $P_d\gtrapprox \frac{m_p}{2}=70$ for $\alpha_P^{[m],M}$.

\subsection{Collocation approach}\label{sec:col2}
Contrary to the Galerkin approach, where the coefficients, evaluated by formula \eqref{eq:ci}, can be computed independently one of each other, in the collocation approach the number of degree of freedom, \ie the $P+1$ coefficients, should match the number of interpolation points. Hence the interpolant of total degree $P_d$, $\alpha_P^{[t]}$, is computed on Padua points with $m_p=P_d$, while the interpolant of maximal degree $P_d$,  $\alpha_P^{[m]}$, is evaluated on tensor products Chebyshev grid with $m_c=P_d$. The matching between interpolant polynomial and interpolating points is strengthened by the following theorem, which provides the near-best optimal approximation in $L^\infty$  associated to a multivariate polynomial degree.

\begin{theorem}\label{th:nearbestd2}
	Let ${\alpha_P^{[t]}}^\star$, and $\alpha_P^{[t]}$,  be the best polynomial approximation of total degree less than or equal to $P_d$ and the polynomial interpolant on $m_p=P_d$ Padua points, respectively. Then 
	\begin{equation}
	\|\alpha - \alpha_P^{[t]}\|_\infty\leq\left(1+\mathcal{O}(\log^2(P_d))\right) \left\|\alpha -{\alpha_P^{[t]}}^\star\right\|_\infty.
	\end{equation}
Analogously, if  ${\alpha_P^{[m]}}^\star$, and $\alpha_P^{[m]}$,  are the best polynomial interpolant of maximal degree less than or equal to $P_d$ and the polynomial interpolant on $m_c=P_d$ tensor product Chebyshev grid, respectively, then 
	\begin{equation}
	\|\alpha - \alpha_P^{[m]}\|_\infty\leq\left(1+\mathcal{O}(\log^2(P_d))\right) \left\|\alpha -{\alpha_P^{[m]}}^\star\right\|_\infty.
	\end{equation}
\end{theorem}
The near-best optimality is derived by the growth of the Lebesgue constant, given in \cite{Bos2006}. Moreover, the following theorem, stated in \cite{Caliari2008}, provides a generalization of the error bounds of Theorem ~\ref{th:convcollocation}, for total degree polynomial interpolant evaluated in Padua points.

\begin{theorem}  
	\label{th:convcoll2}
	If $\alpha$ is continuous and differentiable up to the $k$th derivative in $[-1,1]^2$, with $0< k < \infty$, then the interpolant $\alpha_P^{[t]}$ on $m_p=P_d$ Padua points satisfies the following relation
	\begin{equation}
		\|\alpha - \alpha_P^{[t]}\|_\infty\leq\mathcal{O}\left(\frac{\log^2(P_d))}{P_d^k}\right).
	\end{equation}
\end{theorem}  

Theorem~\ref{th:convcoll2} provides an error bounds  only for the \sae{SAE} case, where it ensures a convergence faster than $\mathcal{O}(P^{-k})$ for all $k\in\N$. 

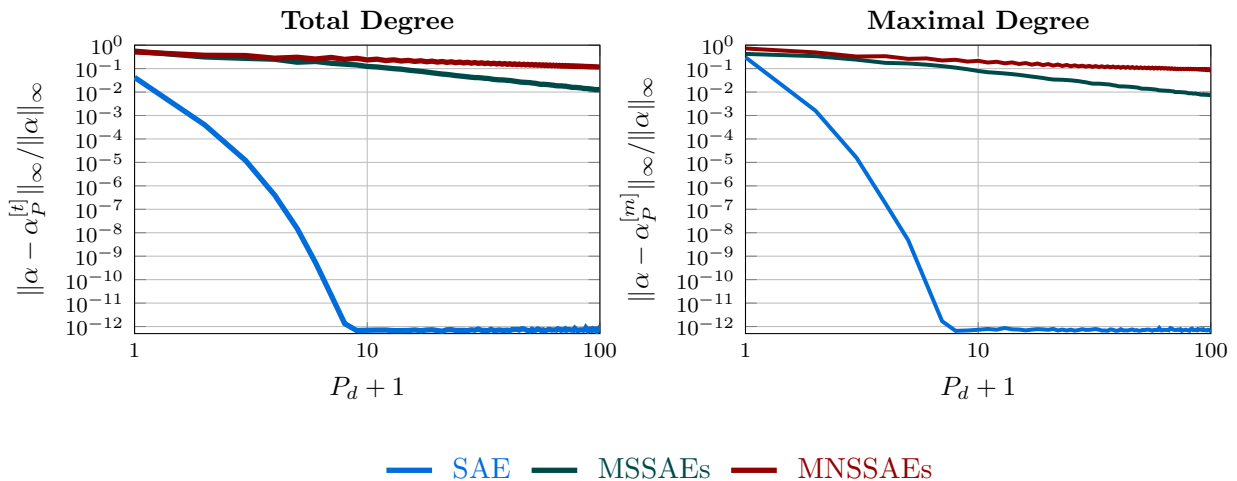
\begin{figure}[!h]
	\begin{center}
		\begin{subfigure}[h]{0.45\linewidth}
			\centering
			\begin{tikzpicture}

\begin{axis}[%
width=0.8\textwidth,
height=0.5\textwidth,
at={(0\textwidth,0\textwidth)},
scale only axis,
xmode=log,
xmin=1,
xmax=100,
xtick={  1,  10, 100},
xticklabels={1,10,100},
xminorticks=false,
xlabel={$P_d+1$},
xmajorgrids,
ymode=log,
ymin=5e-13,
ymax=1,
ytick={1e-12,1e-11,1e-10,1e-09,1e-08,1e-07,1e-06,1e-05,1e-04,1e-03,1e-02,1e-01,1},
ylabel={$\|\alpha-\alpha_P^{[t]}\|_\infty/\|\alpha\|_\infty$},
yminorticks=false,
ymajorgrids,
axis background/.style={fill=white},
ylabel style = {font=\small},
xlabel style = {font=\small},
title style={font=\bfseries, yshift=-0.2cm},
title={\small{Total Degree}}
]
\addplot [color=mycolor1,solid,line width=2.0pt,forget plot]
  table[row sep=crcr]{%
1	0.0423079430397621\\
2	0.000402091248519132\\
3	1.20740617793173e-05\\
4	4.04717588162861e-07\\
5	1.50867624654645e-08\\
6	5.22703300110098e-10\\
7	2.17490486842036e-11\\
8	1.33405872226318e-12\\
9	6.7736618805888e-13\\
10	6.93236257635484e-13\\
11	7.06430913183848e-13\\
12	7.08255059112194e-13\\
13	6.66421312488808e-13\\
14	6.69400750838439e-13\\
15	6.77062163737489e-13\\
16	7.1238978988311e-13\\
17	6.63745898460569e-13\\
18	7.10687253683321e-13\\
19	7.30266419980893e-13\\
20	6.75542042130535e-13\\
21	6.29573564736234e-13\\
22	7.30023200523781e-13\\
23	7.37441393965718e-13\\
24	6.84784381500818e-13\\
25	6.82595406386803e-13\\
26	6.84297942586592e-13\\
27	6.70434433531168e-13\\
28	7.23821104367407e-13\\
29	6.96945354356453e-13\\
30	7.81524920567396e-13\\
31	6.86730137157719e-13\\
32	6.8502760095793e-13\\
33	6.64961995746132e-13\\
34	7.61945754269823e-13\\
35	6.8137930910124e-13\\
36	6.9037842901441e-13\\
37	7.16281301196913e-13\\
38	6.86122088514937e-13\\
39	6.86730137157719e-13\\
40	7.03755499155609e-13\\
41	6.6569165411747e-13\\
42	7.51548122478255e-13\\
43	7.43157051207867e-13\\
44	7.59574364562974e-13\\
45	7.07890229926525e-13\\
46	7.97820624193947e-13\\
47	7.99462355529457e-13\\
48	7.2138890979628e-13\\
49	6.74569164302084e-13\\
50	6.55962875832962e-13\\
51	7.59087925648749e-13\\
52	6.78095846430218e-13\\
53	7.01931353227263e-13\\
54	7.86328504845371e-13\\
55	7.38535881522725e-13\\
56	6.70556043259725e-13\\
57	7.21145690339167e-13\\
58	6.69218336245605e-13\\
59	6.94269940328213e-13\\
60	7.41454515008078e-13\\
61	7.47413391707339e-13\\
62	7.05336425626841e-13\\
63	7.15430033097018e-13\\
64	7.17132569296807e-13\\
65	7.04120328341278e-13\\
66	7.97273380415443e-13\\
67	7.58966315920192e-13\\
68	6.92932233314093e-13\\
69	7.17497398482476e-13\\
70	7.12633009340222e-13\\
71	7.90220016159175e-13\\
72	7.0302584078427e-13\\
73	7.51548122478255e-13\\
74	7.42913831750754e-13\\
75	7.11781741240328e-13\\
76	7.05214815898285e-13\\
77	6.9803984191346e-13\\
78	8.30472836311327e-13\\
79	7.42427392836528e-13\\
80	6.76150090773317e-13\\
81	6.82595406386803e-13\\
82	7.17619008211033e-13\\
83	7.24550762738745e-13\\
84	7.11173692597546e-13\\
85	6.96215695985115e-13\\
86	7.1871349576804e-13\\
87	8.7704936234841e-13\\
88	7.58479877005967e-13\\
89	6.60219216332435e-13\\
90	7.17497398482476e-13\\
91	7.5203456139248e-13\\
92	7.08984717483532e-13\\
93	7.81828944888786e-13\\
94	7.21632129253392e-13\\
95	7.51426512749699e-13\\
96	7.33063443737689e-13\\
97	6.413697084062e-13\\
98	7.28442274052548e-13\\
99	8.15028400784671e-13\\
100	7.29658371338112e-13\\
};
\addplot [color=mycolor2,solid,line width=2.0pt,forget plot]
  table[row sep=crcr]{%
1	0.554834136686369\\
2	0.313799172491642\\
3	0.273805592291579\\
4	0.250360036333426\\
5	0.187102519346337\\
6	0.201390760133319\\
7	0.168994911475212\\
8	0.154188406132532\\
9	0.143375641903148\\
10	0.123781528400162\\
11	0.117414064940436\\
12	0.10766135090475\\
13	0.0994726544911208\\
14	0.0944436215974105\\
15	0.086125749458244\\
16	0.0804615508814851\\
17	0.0764789502830114\\
18	0.0712681542473355\\
19	0.0663692251910699\\
20	0.0621132457228595\\
21	0.0586107110046402\\
22	0.0563131297726932\\
23	0.0534973119084365\\
24	0.0503705646953278\\
25	0.0484176586995323\\
26	0.0464228913256606\\
27	0.04451997359459\\
28	0.0428533239153802\\
29	0.0413936016039186\\
30	0.0412767014816499\\
31	0.0401829598426256\\
32	0.0385850942274417\\
33	0.0380635098940327\\
34	0.0362743658687627\\
35	0.034519466584075\\
36	0.0337631511108624\\
37	0.0329968927917222\\
38	0.033032556765272\\
39	0.0324171609016521\\
40	0.0314668703728051\\
41	0.0310958384470203\\
42	0.03012316591299\\
43	0.0289083396224131\\
44	0.0280920580536268\\
45	0.0271778667361207\\
46	0.0263929535201505\\
47	0.0261734946745691\\
48	0.0256987847155001\\
49	0.0255430236487805\\
50	0.0251914258309773\\
51	0.0244256772372094\\
52	0.0238904409754429\\
53	0.0234303656864434\\
54	0.0225248679487324\\
55	0.0217077328583647\\
56	0.0212942542415111\\
57	0.0214710029556594\\
58	0.0212079623827955\\
59	0.0209494366689653\\
60	0.020751848894667\\
61	0.0202541889179237\\
62	0.0198507451676182\\
63	0.0196012709038447\\
64	0.0190161917061338\\
65	0.0185541126135525\\
66	0.0183976474658359\\
67	0.0183212385882082\\
68	0.0183774407132111\\
69	0.0181374354187246\\
70	0.0178547045483837\\
71	0.0177856539792655\\
72	0.0174024517240522\\
73	0.0169176451158165\\
74	0.0166631940305547\\
75	0.0161627944385495\\
76	0.0158883700767044\\
77	0.0158657830696925\\
78	0.0156931057759775\\
79	0.0157164750270019\\
80	0.0155262605862798\\
81	0.0153245620872968\\
82	0.0152310016375014\\
83	0.0149893118253774\\
84	0.0146261172929124\\
85	0.0143248074951606\\
86	0.0141756582317992\\
87	0.0141612373471792\\
88	0.0139853368750994\\
89	0.01386853593944\\
90	0.0137435299872637\\
91	0.0136240960445789\\
92	0.0134226428118812\\
93	0.0133053916303384\\
94	0.0130787876021784\\
95	0.012743615202645\\
96	0.0126773035074546\\
97	0.0127080763598289\\
98	0.0126752044413397\\
99	0.012648788593362\\
100	0.0124655822983376\\
};
\addplot [color=mycolor3,solid,line width=2.0pt,forget plot]
  table[row sep=crcr]{%
1	0.514085179727774\\
2	0.369982771551656\\
3	0.360697885163594\\
4	0.277806230313123\\
5	0.304532885481923\\
6	0.25571821978925\\
7	0.299223376752894\\
8	0.247071623957016\\
9	0.275009679091159\\
10	0.23433475475605\\
11	0.252334518296185\\
12	0.220799636283729\\
13	0.236053115574438\\
14	0.211056662698057\\
15	0.224168852805607\\
16	0.203115796056398\\
17	0.214896672227846\\
18	0.195688612198107\\
19	0.207040331482016\\
20	0.189365866442901\\
21	0.200597348682871\\
22	0.183870996817613\\
23	0.195025005444229\\
24	0.179119953525703\\
25	0.189592054478702\\
26	0.174644053073998\\
27	0.185055848816117\\
28	0.170675745881102\\
29	0.1807468897969\\
30	0.167327173230903\\
31	0.176514911603662\\
32	0.163821950280826\\
33	0.172440734740452\\
34	0.160880637293071\\
35	0.169284083969824\\
36	0.158144950579737\\
37	0.166371133574865\\
38	0.155271878563392\\
39	0.163721149282215\\
40	0.152747187192704\\
41	0.161048938533065\\
42	0.15030832664843\\
43	0.158457418585816\\
44	0.147958468409793\\
45	0.155888448254538\\
46	0.146089522199141\\
47	0.15333508659073\\
48	0.144276630909863\\
49	0.150840121604186\\
50	0.142501031231271\\
51	0.14889432690228\\
52	0.140769996864609\\
53	0.147178147748562\\
54	0.139139494949647\\
55	0.145510907949837\\
56	0.137552892985139\\
57	0.144022479545699\\
58	0.136016713372889\\
59	0.14254123621413\\
60	0.134547036877318\\
61	0.14106828836529\\
62	0.133161784706594\\
63	0.13961945703175\\
64	0.131776148864552\\
65	0.138181589881963\\
66	0.130439211409004\\
67	0.136731611375868\\
68	0.129208029981176\\
69	0.1351541333963\\
70	0.128091475821689\\
71	0.133893429765045\\
72	0.12708501919249\\
73	0.132637698384495\\
74	0.126106669940285\\
75	0.1313896766694\\
76	0.125157928185848\\
77	0.130160978451293\\
78	0.124235650450377\\
79	0.128949608476929\\
80	0.123318085404942\\
81	0.127750228898856\\
82	0.122430563106753\\
83	0.126586316083949\\
84	0.121587501309744\\
85	0.125435120764963\\
86	0.12077152610356\\
87	0.124306406465933\\
88	0.119987242043642\\
89	0.123187383086031\\
90	0.119231029128801\\
91	0.122083307654032\\
92	0.118497081467575\\
93	0.120995293205194\\
94	0.117805299149732\\
95	0.119924780020168\\
96	0.117111611433399\\
97	0.118872432718318\\
98	0.116416990267313\\
99	0.11783815977769\\
100	0.115736537039547\\
};
\end{axis}
\end{tikzpicture}%
		\end{subfigure}\quad
		\begin{subfigure}[h]{0.45\linewidth}
			\centering
			\begin{tikzpicture}

\begin{axis}[%
width=0.8\textwidth,
height=0.5\textwidth,
at={(0\textwidth,0\textwidth)},
scale only axis,
xmode=log,
xmin=1,
xmax=100,
xtick={  1,  10, 100},
xticklabels={1,10,100},
xminorticks=false,
xlabel={$P_d+1$},
xmajorgrids,
ymode=log,
ymin=5e-13,
ymax=1,
yminorticks=false,
ytick={1e-12,1e-11,1e-10,1e-09,1e-08,1e-07,1e-06,1e-05,1e-04,1e-03,1e-02,1e-01,1},
ylabel={$\|\alpha-\alpha_P^{[m]}\|_\infty/\|\alpha\|_\infty$},
ymajorgrids,
axis background/.style={fill=white},
ylabel style = {font=\small},
xlabel style = {font=\small},
title style={font=\bfseries, yshift=-0.2cm},
title={\small{Maximal Degree}}
]
\addplot [color=mycolor1,solid,line width=1.5pt,forget plot]
  table[row sep=crcr]{%
1	0.303583574173491\\
2	0.00159720604220022\\
3	1.58023010888742e-05\\
4	1.72581822909278e-07\\
5	4.86777731090139e-09\\
6	6.48786685750852e-11\\
7	1.68891591019062e-12\\
8	6.53652290990392e-13\\
9	6.79555163172894e-13\\
10	7.27834225409766e-13\\
11	7.90463235616288e-13\\
12	7.31968956180682e-13\\
13	8.67442193792458e-13\\
14	7.60912071577094e-13\\
15	7.43521880393536e-13\\
16	6.86608527429163e-13\\
17	7.39022320436951e-13\\
18	7.76842946017976e-13\\
19	7.56594926213343e-13\\
20	7.23821104367407e-13\\
21	7.92895430187414e-13\\
22	6.8071045559418e-13\\
23	7.36954955051493e-13\\
24	6.79676772901451e-13\\
25	7.3172573672357e-13\\
26	6.6970477515983e-13\\
27	7.07586205605134e-13\\
28	7.69424752576039e-13\\
29	7.19321544410822e-13\\
30	7.21996958439061e-13\\
31	7.50940073835473e-13\\
32	7.81828944888786e-13\\
33	7.37319784237162e-13\\
34	6.75055603216309e-13\\
35	6.76879749144655e-13\\
36	6.58516680132646e-13\\
37	7.54223536506495e-13\\
38	7.23577884910294e-13\\
39	6.61070484432329e-13\\
40	6.69461555702717e-13\\
41	7.69303142847482e-13\\
42	6.81014479915571e-13\\
43	6.86486917700607e-13\\
44	6.77852626973106e-13\\
45	7.20780861153498e-13\\
46	7.47170172250226e-13\\
47	6.85878869057825e-13\\
48	7.0339066996994e-13\\
49	6.60340826060991e-13\\
50	6.74447554573528e-13\\
51	7.15673252554131e-13\\
52	7.0570125481251e-13\\
53	7.33306663194802e-13\\
54	7.4005600312968e-13\\
55	7.06430913183848e-13\\
56	6.65266020067523e-13\\
57	7.11781741240328e-13\\
58	7.16281301196913e-13\\
59	7.07160571555186e-13\\
60	8.65435633271278e-13\\
61	6.9037842901441e-13\\
62	6.89527160914515e-13\\
63	7.23942714095963e-13\\
64	7.08741498026419e-13\\
65	7.19929593053603e-13\\
66	6.91959355485642e-13\\
67	8.58443073879288e-13\\
68	8.16487717527347e-13\\
69	7.30144810252337e-13\\
70	7.11295302326102e-13\\
71	7.0570125481251e-13\\
72	7.44129929036317e-13\\
73	7.19807983325047e-13\\
74	7.67965435833362e-13\\
75	6.57422192575638e-13\\
76	6.90621648471523e-13\\
77	7.11416912054659e-13\\
78	7.34279541023253e-13\\
79	7.01931353227263e-13\\
80	6.82352186929691e-13\\
81	6.62286581717893e-13\\
82	7.76842946017976e-13\\
83	8.12474596484987e-13\\
84	7.08741498026419e-13\\
85	7.06066083998179e-13\\
86	7.28442274052548e-13\\
87	6.94877988970995e-13\\
88	7.80612847603223e-13\\
89	7.24550762738745e-13\\
90	7.37198174508605e-13\\
91	7.53554682999435e-13\\
92	7.06795742369517e-13\\
93	7.5969597429153e-13\\
94	7.10079205040539e-13\\
95	6.84541162043705e-13\\
96	6.8004160208712e-13\\
97	7.22726616810399e-13\\
98	7.19807983325047e-13\\
99	6.71528921088175e-13\\
100	6.94756379242439e-13\\
};
\addplot [color=mycolor2,solid,line width=1.5pt,forget plot]
  table[row sep=crcr]{%
1	0.422550656218589\\
2	0.339250377723348\\
3	0.240577058266873\\
4	0.175241590559657\\
5	0.16564060609695\\
6	0.146878498881937\\
7	0.127751324571409\\
8	0.109910341003431\\
9	0.093591934054632\\
10	0.0791845478859808\\
11	0.0701335502730012\\
12	0.0662401689463481\\
13	0.0621287426295505\\
14	0.0580236558644545\\
15	0.0539375876418451\\
16	0.0495128509139402\\
17	0.0456954937923777\\
18	0.0426713221077203\\
19	0.0396987255378564\\
20	0.0365003372540195\\
21	0.0342572837599574\\
22	0.033764696688378\\
23	0.0333667656031876\\
24	0.0325035910160959\\
25	0.0316356750277112\\
26	0.0303325684943277\\
27	0.0291405564040028\\
28	0.0273946628514051\\
29	0.0261765434629324\\
30	0.0243279734292826\\
31	0.0228264414713277\\
32	0.0227678256603442\\
33	0.0227072882777788\\
34	0.022354014324996\\
35	0.0217678160392136\\
36	0.0213222493663716\\
37	0.0206823315942495\\
38	0.0195341873325172\\
39	0.0189098155714243\\
40	0.0177503862626037\\
41	0.0173691361908889\\
42	0.0173956404484834\\
43	0.0171836777511115\\
44	0.0169818969643228\\
45	0.0167288090634677\\
46	0.0163805674707778\\
47	0.0160165390422777\\
48	0.0156720030947424\\
49	0.0151952276106296\\
50	0.0141784774480056\\
51	0.0138893375456803\\
52	0.0138695136232669\\
53	0.0137998010227377\\
54	0.0136397870903763\\
55	0.0134409260866743\\
56	0.0132459735521696\\
57	0.0130279013376532\\
58	0.0127614820673092\\
59	0.0123871847214636\\
60	0.0118036068911761\\
61	0.0117723650398068\\
62	0.0115474073559229\\
63	0.0115553374983689\\
64	0.0113929682033329\\
65	0.0112634673135713\\
66	0.0112735611772334\\
67	0.0110360888567221\\
68	0.0107658826327087\\
69	0.0107160200651655\\
70	0.0100767108998616\\
71	0.0100410758484651\\
72	0.00999440750769083\\
73	0.00995585768314431\\
74	0.00996127724049469\\
75	0.00977526989661378\\
76	0.00973482834829514\\
77	0.00966022905708747\\
78	0.00953217764184804\\
79	0.00934670502663198\\
80	0.00888125140556434\\
81	0.00887610357404503\\
82	0.00888948527224124\\
83	0.00883409398244294\\
84	0.00856056196329614\\
85	0.00862448332006309\\
86	0.00875525004127076\\
87	0.00840615127182865\\
88	0.00828338156093756\\
89	0.00824365699295192\\
90	0.00784829729412411\\
91	0.00773876175089132\\
92	0.00771611714773281\\
93	0.0079891848378423\\
94	0.0076697591052901\\
95	0.00762196691925062\\
96	0.00780059152991016\\
97	0.00750877832772247\\
98	0.00748910100002668\\
99	0.0074028806306977\\
100	0.00720071101563285\\
};
\addplot [color=mycolor3,solid,line width=1.5pt,forget plot]
  table[row sep=crcr]{%
1	0.722259166329192\\
2	0.490422648716763\\
3	0.327440711651819\\
4	0.337602464323005\\
5	0.260078851964008\\
6	0.27180276947275\\
7	0.222024025669294\\
8	0.23787655295202\\
9	0.199262884641781\\
10	0.214048206696258\\
11	0.18170459909006\\
12	0.195960020919516\\
13	0.168181810293638\\
14	0.183173348708841\\
15	0.157222365395695\\
16	0.172126674204142\\
17	0.148025108208493\\
18	0.163828766022945\\
19	0.140484624055032\\
20	0.155919863439461\\
21	0.133908236573973\\
22	0.148843391722011\\
23	0.129551531556528\\
24	0.142611231681209\\
25	0.126268529377106\\
26	0.137119678770248\\
27	0.124644821002278\\
28	0.132396455235091\\
29	0.12157567279223\\
30	0.129836105543407\\
31	0.118493107171059\\
32	0.127495404174405\\
33	0.117204025463057\\
34	0.125142153995739\\
35	0.11526612507914\\
36	0.124139756709221\\
37	0.113590454579382\\
38	0.121915320021521\\
39	0.111286949114859\\
40	0.1191043074826\\
41	0.108582469915742\\
42	0.118537047288689\\
43	0.108760970594473\\
44	0.117760353766705\\
45	0.106824500036776\\
46	0.116284733451648\\
47	0.105323486566297\\
48	0.113895283925833\\
49	0.105435164867205\\
50	0.111896507527187\\
51	0.104576808203868\\
52	0.110572716726642\\
53	0.103445874501047\\
54	0.110661492287757\\
55	0.101556280661402\\
56	0.110012637676153\\
57	0.100781548451117\\
58	0.108611962631958\\
59	0.100878854677121\\
60	0.107956915865676\\
61	0.100505192128891\\
62	0.106878767118497\\
63	0.100269428360305\\
64	0.105206379711037\\
65	0.0989247963721322\\
66	0.10352476504157\\
67	0.09764259866072\\
68	0.10237505001075\\
69	0.0963930994574902\\
70	0.100748115662715\\
71	0.0949481448013955\\
72	0.099102727899164\\
73	0.0946887323924571\\
74	0.0987483875258705\\
75	0.0949016650665085\\
76	0.0981226521039142\\
77	0.0951634601673914\\
78	0.098426035333183\\
79	0.0948136292167842\\
80	0.0984498051714126\\
81	0.0937802470016091\\
82	0.0981608188580897\\
83	0.0935827756003781\\
84	0.0975709084198368\\
85	0.0933232190027648\\
86	0.0966874722907523\\
87	0.0924954981773549\\
88	0.0956323366113872\\
89	0.0910618953761862\\
90	0.0955975679047428\\
91	0.090102429791742\\
92	0.0953412274025991\\
93	0.0896993695897819\\
94	0.0948772982280845\\
95	0.0888172111878545\\
96	0.094191370359377\\
97	0.087437011919688\\
98	0.0932986734484433\\
99	0.0862645910864327\\
100	0.0922040574382094\\
};
\end{axis}
\end{tikzpicture}%
		\end{subfigure}
		\vspace{-0.01\textwidth}  
		\begin{tikzpicture}
		\hspace{0.05\linewidth}
		\begin{customlegend}[legend columns=5,legend style={align=left,draw=none,column sep=1ex},legend entries={\sae{SAE},\mssae{MSSAEs},\mnssae{MNSSAEs}}]
		\addlegendimage{color=mycolor1, mark=none,solid,line width=2.0pt,line legend}
		\addlegendimage{color=mycolor2, mark=none,solid,line width=2.0pt}   
		\addlegendimage{color=mycolor3, mark=none,solid,line width=2.0pt}
		\end{customlegend}
		\end{tikzpicture}
	\end{center}
	\caption{Convergence rates of interpolant polynomial for collocation approach of the spectral abscissa functions of Example~\ref{ex:d2}, based on Padua points for total degree and tensor product Chebyshev grid for maximal degree.  } 
	\label{fig:collocationd2}
\end{figure}

Figure~\ref{fig:collocationd2} shows the experiments on the spectral abscissa  functions of Example~\ref{ex:d2}, for which we used the software Chebfun relying on Chebyshev expansion. The convergence rates are analogous to the ones of the Galerkin approach (Figure~\ref{fig:gd2_approx_err}) the only difference is that the numerical error is negligible for the collocation approach.

\section{Conclusions} \label{sec:conclusions}
This paper, other than explaining the parallelism between polynomial approximation and PC theory,  analyzes the approximation of polynomial series \eqref{eq:spabsapp} of $\alpha(\omega)$ (and its PC expansion \eqref{eq:spabspce} of $\alpha(\bomega)$) w.r.t.~the behavior of spectral abscissa function \eqref{eq:spabsfun}. The analyses show that the lack of smoothness properties heavily affects both the approximation errors of the Galerkin and collocation approaches, and the numerical errors in the approximation of the coefficients $c_i$ by integration methods.

The convergence rates between the Galerkin and collocation approaches are similar, if the numerical errors are negligible. In particular, for the test-examples analyzed (Examples~\ref{ex:d1} and \ref{ex:d2}), if the spectral abscissa behaves smoothly (\sae{SAE}), then the polynomial approximation convergences with an order of $\mathcal{O}(P_d^{-P_d})$ where $P_d=P$ for $D=1$. However, the non-differentiable and non-Lipschitz continuous behaviors (represented in the benchmark examples by \mssae{MSSAEs} and \mnssae{MNSSAEs}) present in the univariate case a order of convergence of $\mathcal{O}(P^{-1})$ and $\mathcal{O}(P^{-0.5})$, respectively, while in the bivariate case they converge approximately as  $\mathcal{O}(P_d^{-1})$ and $\mathcal{O}(P_d^{-0.3})$, respectively. These latter cases are not deeply studied in the literature on the spectral abscissa approximation, even though they easily occur when applying stability optimization in the context of the design of the controllers (cf., \eg \cite{Fenzi2017a,Fenzi2017} and reference therein).

The present work, moreover, reviews the main theorems on univariate and bivariate polynomial approximation on Chebyshev and Legendre bases for differentiable functions, providing convergence rates for benchmark examples with non-differentiable and non-Lipschitz continuous functions, \ie \mssae{MSSAEs} and \mnssae{MNSSAEs}, respectively.

A last contribution of this paper occurs in the advice of Section~\ref{sec:coupling}, which correlates the  polynomial degree of $\alpha_P^M$ with the number of points $M+1$ of the integration methods, in such a way that the numerical error of the integration method does not affect the quality of the approximation of the Galerkin approach.

\section*{Acknowledgments} 
The authors would like to thank A.~Bultheel for the careful proofreading and his advice, and L.~N.~Trefethen for pointing to valuable references.

This work was supported by the project C14/17/072 of the KU Leuven Research Council, by the project G0A5317N of the Research Foundation-Flanders (FWO - Vlaanderen), and by the project UCoCoS, funded by the European Unions Horizon 2020 research and innovation program under the Marie Sklodowska-Curie Grant Agreement No 675080.

\appendix
\section{Galerkin approach for the spectral abscissa functions in Example~\ref{ex:d1}}\label{app:ex1}    

The spectral abscissa functions for eigenvalue problems \eqref{eq:SAE1}, \eqref{eq:MSSAE1}, and \eqref{eq:MNSSAE1} are given by
\begin{align}
	\begin{split}
		\alpha^{\text{I}}(\omega)=e^\omega,\  \omega\in[-1,1];\quad  
		\alpha^{\text{II}}(\omega)=\begin{cases} 0, & \text{ if }\omega\in[-1,0),\\
			\omega, &\text{ if }\omega\in[0,1];
		\end{cases}\quad
		\alpha^{\text{III}}(\omega)=\begin{cases} 0, &\text{ if }\omega\in[-1,0),\\
			\sqrt{\omega}, &\text{ if }\omega\in[0,1];
		\end{cases}
	\end{split}
	\label{eq:spabsex1}
\end{align}
respectively. In this appendix, we analytically compute the quantities, which appear in  the evaluation of the coefficients $c_i$ in equation  \eqref{eq:ci}.

The $\rho$-norms of the Legendre polynomial are given by $\|p_i\|_{\rho}^2=\frac{1}{2i+1}$ for all $i\in\N$ (see \eg 22.2.10 in \cite{Abramowitz1965}).

To analytically express the coefficients $c_i$ in  \eqref{eq:ci}, the evaluation of the $\rho$-inner product of spectral abscissas \eqref{eq:spabsex1} and the $i$th Legendre polynomials is needed. We start analyzing  $\alpha^{\text{I}}(\omega)$ and we furnish the corresponding $\rho$-inner product iteratively, starting from:
\begin{equation*}
	\langle \alpha^{\text{I}},p_0\rangle_{\rho}=\frac{1}{2}\left(e-e^{-1}\right).
\end{equation*}
Set $i\in\N$ and $i\geq1$, by using integration by parts, and  the following Legendre polynomial properties
\begin{align*}
p_{i}(\pm1)=(\pm1)^{i}, \quad\text{and } \frac{\diff p_{i+1}(\omega)}{\diff\omega}=\sum_{k=0}^{\lfloor i/2\rfloor}\frac{p_{i-2k}(\omega)}{\|p_{i-2k}\|_{\rho}^2},\quad \text{ where } \left\lfloor\frac{i}{2}\right\rfloor=\max_{j\in\N}\left\{j\leq \frac{i}{2}\right\},	
\end{align*}
we get
\begin{align*}
	\langle\alpha^{\text{I}},p_{i+1}\rangle_{\rho}&=\frac{1}{2}\left(e+\frac{(-1)^i}{e}-2\sum_{k=0}^{\lfloor i/2\rfloor}\frac{\langle\alpha^{\text{I}},p_{i-2k}\rangle_{\rho}}{\|p_{i-2k}\|_{\rho}^2}\right)=\frac{1}{2}\left(e+\frac{(-1)^i}{e}-2\sum_{k=0}^{\lfloor i/2\rfloor} {c}_{i}^{\text{I}}\right),
	\end{align*}
where ${c}_{i}^{\text{I}}$ indicates the $i$th coefficients of polynomial approximation \eqref{eq:spabsapp} of $\alpha^{\text{I}}(\omega)$.

For all $i\in\N$, the $\rho$-inner product of the spectral abscissas $\alpha^{\text{II}}(\omega)$ and $\alpha^{\text{III}}(\omega)$ can be evaluated by using the relation 22.13.8 and 22.13.9 in \cite{Abramowitz1965}:
\begin{align*}
	\langle\alpha^{\text{II}},p_{i}\rangle_{\rho}&=\begin{cases}
		\frac{(-1)^j\Gamma\left(j-\frac{1}{2}\right)}{4\Gamma\left(-\frac{1}{2}\right)\Gamma\left(j+2\right)}, &{\text{if }i=2j},\\
		\frac{1}{6}, &\text{if }i=1,\\
		0, &\text{if }i=2j+1,\ i>1
	\end{cases}\quad
	\langle\alpha^{\text{III}},p_{i}\rangle_{\rho}&=\begin{cases}
		\frac{(-1)^j\Gamma\left(j-\frac{1}{4}\right)\Gamma\left(\frac{3}{4}\right)}{4\Gamma\left(-\frac{1}{4}\right)\Gamma\left(j+\frac{7}{4}\right)}, &\text{if }i=2j,\\
		\frac{(-1)^j\Gamma\left(j+\frac{1}{4}\right)\Gamma\left(\frac{5}{4}\right)}
		{4\Gamma\left(\frac{1}{4}\right)\Gamma\left(j+\frac{9}{4}\right)}, &\text{if }i=2j+1,
	\end{cases}
\end{align*}
where $\Gamma(\cdot)$ denotes the Gamma function. The two last inner product are, in fact, rational numbers which can be computed without any error via a symbolic software.

\end{document}